\documentclass [a4paper,twoside,11pt]{article}
\usepackage[utf8]{inputenc}
\usepackage{amsmath,amsfonts,amssymb}
\usepackage{vmargin,graphicx,theorem}
\usepackage[english]{babel}
\usepackage{enumerate}
\RequirePackage[colorlinks,linkcolor=blue,citecolor=red,urlcolor=blue]{hyperref}
\usepackage{color}

\setpapersize[portrait]{A4}
\setmarginsrb{1.5cm}{1cm}{1.5cm}{1.5cm}{1.5cm}{0cm}{0.5cm}{2cm}
\usepackage{pst-fill,pst-grad,pst-plot,pst-eucl,pstricks-add,pst-node}

\newcommand{\rr}{\ensuremath{\mathbb{R}}}
\newcommand{\dR}{\ensuremath{\mathbb{R}}}
\newcommand{\N}{\ensuremath{\mathbb{N}}}
\usepackage{ulem}
\newcommand{\R}{\dR}

\newtheorem{thm}{Theorem}[section]
\newtheorem{coro}[thm]{Corollary}
\newtheorem{prop}[thm]{Proposition}
\newtheorem{lemma}[thm]{Lemma}

\newtheorem{remark}[thm]{Remark}

\newcommand{\proofend}{~$\rhd$}
\newcommand{\proofbegin}{~$\lhd$}

\newenvironment{eproof}
               {\noindent {\it{\textbf{Proof}}}\\\proofbegin~}
               {\proofend\\}

\newcommand{\PAR}[1]{\ensuremath{{\left(#1\right)}}} 
\newcommand{\SBRA}[1]{\ensuremath{{\left[#1\right]}}} 
\newcommand{\BRA}[1]{\ensuremath{{\left\{#1\right\}}}} 
\newcommand{\entf}[1]{{\rm{Ent}}_{#1}}
\newcommand{\varf}[1]{{\rm{Var}}_{#1}}
\newcommand{\He}[1]{{\rm{Hess}}({#1})}

\newcommand{\e}{\varepsilon}

\newcommand{\tr}{\textrm{tr}}

\newcommand{\beq}{\begin{equation}}\newcommand{\eeq}{\end{equation}}
\begin{document}

\title{Dimensional improvements of the logarithmic Sobolev, Talagrand and Brascamp-Lieb inequalities}

\author{ Fran{\c c}ois Bolley\thanks{Laboratoire de Probabilit\'es et Mod\`eles Al\'eatoires, Umr Cnrs 7599, Universit\'e Pierre et Marie Curie, Paris, France. francois.bolley@upmc.fr}, Ivan Gentil\thanks{Univ Lyon, Universit\'e Claude Bernard Lyon 1, CNRS UMR 5208, Institut Camille Jordan, 43 blvd. du 11 novembre 1918, F-69622 Villeurbanne cedex, France. gentil@math.univ-lyon1.fr}\, and Arnaud Guillin\thanks{Laboratoire de Math\'ematiques, Umr Cnrs 6620, Universit\'e Blaise Pascal, Clermont-Ferrand, France. guillin@math.univ-bpclermont.fr}}

\date{\today}

\maketitle

\abstract{\noindent In this work we consider dimensional improvements of the logarithmic Sobolev, Talagrand and Brascamp-Lieb inequalities. For this we use optimal transport methods and the Borell-Brascamp-Lieb inequality. These refinements can be written as a deficit in the classical inequalities. They have the right scale with respect to the dimension. They lead to sharpened concentration properties as well as refined contraction bounds, convergence to equilibrium and short time behavior for the laws of solutions to stochastic differential equations.}

\bigskip

\bigskip


\bigskip


\section*{Introduction}
\parindent=0pt

We shall be concerned with diverse ways of measuring and bounding the distance between probability measures, and the links between them. We will focus on three main inequalities that we now describe. 

\begin{itemize}

\item A probability measure $\mu$ on $\R^n$ satisfies a logarithmic Sobolev inequality (in short LSI) with constant $R>0$ (see~\cite{bgl-book} for instance) if for all probability measures $\nu$ in $\R^n$,  absolutely continuous with respect to $\mu$,
\begin{equation}\label{LS}
H(\nu|\mu) \le \frac1{2R}\, I(\nu|\mu).
\end{equation}
Here $H$ and $I$ are the relative entropy and the Fisher information, defined for  $f=\frac{d\nu}{d\mu}$ by
\begin{equation}
\label{eq-def-fisher}
H(\nu|\mu)= \entf{\mu}(f) = \int f \,  \log f \, d\mu \qquad \textrm{and} \qquad 
I(\nu|\mu)=\int \frac{|\nabla f|^2}fd\mu.
\end{equation}
For $I$ we assume that $\nabla f / f \in L^2(\nu).$
\item A probability measure $\mu$ in $\R^n$ satisfies a Talagrand transportation inequality~\cite{talagrand96} with constant $R>0$ if for all $\nu$ absolutely continuous with respect to $\mu$
\begin{equation}\label{talaclass}
 W_2^2(\nu,\mu)\le \frac{2}{R}H(\nu|\mu).
\end{equation}
Here $W_2$ is  the Monge-Kantorovich-Wasserstein distance; it is defined for $\mu$ and $\nu$ in $P_2(\R^n)$ by
$$
W_2(\mu, \nu) = \inf_{\pi} \left( \iint \vert y-x \vert^2 \, d\pi(x,y) \right)^{1/2}
$$
where $\pi$ runs over the set of (coupling) measures on $\rr^{n} \times \rr^n$ with respective marginals $\mu$ and $\nu.$ 
We let $P_2(\R^n)$ be the space of probability measures $\mu$ on $\rr^n$ with finite second moment, that is, 
$\int \vert x \vert^2 d\mu(x)<+\infty$ (see~\cite{ambrosio-gigli-savare},~\cite{villani-book1}).

By the Otto-Villani Theorem~\cite{ov00}, the logarithmic  Sobolev inequality~\eqref{LS} implies the Talagrand inequality~\eqref{talaclass} with the same constant (see also~\cite{bgl},~\cite[Chap.~22]{villani-book1}).
\item Let $\mu$ be a probability measure in $\R^n$ with density $e^{-V}$ where $V$ is a $\mathcal C^2$ and strictly convex function. Then the Brascamp-Lieb inequality asserts that for all smooth functions $f$, 
\begin{equation}\label{BL}
\varf{\mu}(f) \leq \int \nabla f \cdot \He{V}^{-1} \nabla f \, d\mu.
\end{equation}
Here  $\varf{\mu} (f) = \int f^2 d\mu - (\int fd\mu)^2$ is the variance of $f$ under the measure $\mu$, 
see~\cite[Sect~4.9.1]{bgl-book} for instance.

\end{itemize}

\bigskip

The standard Gaussian measure $\gamma$ in $\R^n$ with density $e^{-V}$ for $V(x) = \vert x \vert^2/2 + n \log(2 \pi) /2,$ satisfies the three inequalities~\eqref{LS},~\eqref{talaclass}  with $R=1$ and~\eqref{BL}.  
In fact, in the Gaussian case, the Brascamp-Lieb inequality~\eqref{BL} can be obtained from~\eqref{LS} by linearization, namely by taking $\nu = f \mu$ with $f$ close to~$1$. Let us note that in this case $\He V={\rm Id}_n$, the Brascamp-Lieb inequality becomes exactly the Poincar\'e inequality.  Moreover these inequalities are optimal for the Gaussian measure: by direct computation, equality holds in~\eqref{LS} and~\eqref{talaclass} for translations of $\gamma,$ that is, for measures $\nu=\exp(a \cdot x-\frac{\vert a \vert^2}2) \gamma$ with $a\in\R^n$; equality holds in~\eqref{BL} for $f(x) = b \cdot x$, $b\in\R^n$ (see~\cite[Chap.~4 and~5]{bgl-book}).

\bigskip

Inequalities~\eqref{LS},~\eqref{talaclass} and~\eqref{BL} share the significant property of tensorisation, leading to possible constants $R$ independent of the dimension of the space. In other words, if a probability measure $\mu$ satisfies one of these three inequalities with constant 
$R>0$, then for any  $N\in\N^*$, the product measure $\mu^N=\otimes^N\mu$ satisfies the same inequality with the same constant $R$.  This can be interesting in applications to problems set in large or infinite dimensions. 


However, for regularity or integrability arguments, one may need more precise forms capturing the precise dependence on the dimension. 
Such dimension dependent improvements have been observed in the Gaussian case. Namely, the dimensional improvement
\begin{equation}\label{eq-logsob-bakry-ledoux}
H(\nu \vert \gamma) \leq
\frac{1}{2} \int \vert x \vert^2d\nu - \frac{n}{2} + \frac{n}{2} \log \Big( 1 + \frac{1}{n} \Big(I(\nu \vert \gamma) + n -  \int  \vert x \vert^2d\nu\Big) \Big)
\end{equation}
of the logarithmic  Sobolev inequality~\eqref{LS} has been obtained by D. Bakry and M. Ledoux~\cite{bakryledoux-liyau} by self-improvement from the Euclidean logarithmic  Sobolev inequality, or by semigroup arguments on the Euclidean heat semigroup (see also~\cite[Sect.~6.7.1]{bgl-book} and the early work~\cite{C} by E. Carlen). The dimensional improvement\begin{equation}\label{taladimgauss}
W_2^2(\nu, \gamma) \leq \int \vert x \vert^2d\nu + n - 2 n \exp \PAR{\int \frac{\vert x \vert^2}{2n}d\gamma - \frac12 - \frac{1}{n} H(\nu \vert \gamma) }
\end{equation}
of the Talagrand inequality~\eqref{talaclass} has been derived in~\cite{BBG12}; the argument is based on local hypercontractivity techniques on an associated Hamilton-Jacobi semigroup and fine properties of the heat semigroup. It has further been observed in~\cite{bakryledoux-liyau} that linearizing~\eqref{eq-logsob-bakry-ledoux} leads to the dimensional improvement
 \begin{equation}\label{Pdimgauss}
\varf{\gamma}(f) \leq \int |\nabla f |^2 \, d\gamma- \frac{1}{2n}\Big(\int (\vert x \vert^2-n) fd\gamma \Big)^2
\end{equation}
of the Brascamp-Lieb (or Poincar\'e) inequality~\eqref{BL} for the Gaussian measure (see also \cite[Sect.~6.7.1]{bgl-book}).  On the other hand, by a spectral analysis of the Ornstein-Uhlenbeck semigroup, the bound
\begin{equation}\label{NPdimgauss}
\varf{\gamma}(f) \leq \frac12\int |\nabla f |^2 \, d\gamma+\frac12 \, \Big|\int \nabla fd\gamma\Big|^2
\end{equation}
has been established in~\cite[Sect.~6.2]{gnp}. By the Cauchy-Schwarz inequality, it improves upon~\eqref{BL}. Naturally, both inequalities~\eqref{Pdimgauss} and~\eqref{NPdimgauss} are optimal, and equality holds for $f(x) = a \cdot x$; equality also holds for $f(x) = \vert x \vert^2$, in fact for the first two Hermite polynomials. The above proofs of~\eqref{eq-logsob-bakry-ledoux},~\eqref{taladimgauss} and~\eqref{NPdimgauss} are very specific to the Gaussian case and can not be extended to other measures.

\bigskip 

These dimensional improvements can also be written as a deficit in the classical non dimensional versions~\eqref{LS},~\eqref{talaclass},~\eqref{BL} of the inequalities: namely, for the logarithmic Sobolev  ($LSI$ in short) and Talagrand ($Tal$ in short) inequalities, lower bounds on the quantities 
$$
\delta_{LSI}(\nu|\mu):=\frac12 I(\nu|\mu)- R\, H(\nu|\mu)
\qquad \mathrm{and} \qquad 
\delta_{Tal}(\nu|\mu):=H(\nu|\mu)-\frac{R}2\, W_2^2(\nu,\mu).
$$

The problem of dimensional refinements of standard functional inequalities has been recently considered in an intensive manner. Via the development of refined optimal transportation tools, beautiful results for the Gaussian isoperimetric inequality were obtained by  Figalli-Maggi-Pratelli \cite{FMP10} (see also R. Eldan \cite{E} or \cite{FI13} for convex cones). Further recent results have been established on deficit in the logarithmic Sobolev inequality in the Gaussian case by Figalli-Maggi-Pratelli \cite{FMP13}, Indrei-Marcon \cite{IM14} and Bobkov $\&$ al \cite{bgrs14}. In particular \cite{bgrs14} rediscovers  \eqref{eq-logsob-bakry-ledoux} and extends earlier results obtained in dimension one by  Barthe-Kolesnikov \cite{BK08} on the Talagrand deficit. Fathi-Indrei-Ledoux \cite{FIL14}  also considers these deficits, particularly emphasizing the case where $\nu$ has additional properties, such as a Poincar\'e inequality ensuring a better constant in the logarithmic Sobolev inequality. Very recently D. Cordero-Erausquin \cite{cordero15} has studied refinements of the Talagrand and Brascamp-Lieb inequalities via optimal transport tools.

Let us also quote C. Villani \cite[p.~605]{villani-book1}:
\smallskip

\begin{center}\begin{minipage}{14cm}{\it  There is no well-identified analog of Talagrand inequalities that would take advantage of the finiteness of the dimension to provide sharper concentration inequalities}\end{minipage}
\end{center}

\smallskip

as a motivation to investigate further the problem. As we will see there are other striking applications of these dimensional refinements than sole concentration.

Finally recall that the so-called Bakry-\'Emery criterion (or $\Gamma_2$-criterion) ensures that the measure $\mu$ with density $e^{-V}$ satisfies the logarithmic Sobolev inequality~\eqref{LS}  and Talagrand inequality~\eqref{talaclass} as soon as the potential $V$ satisfies $\He{V}\geq R\,{\rm Id}_n$ with $R>0,$  as symmetric matrices.  One of the goals of this paper is to extend the above dimensional inequalities under this condition with $R>0$ or only $\He{V}>0$. For this we shall use multiple tools and we will compare our inequalities with other recent extensions. Applications to concentration inequalities and short and long time behaviour for the laws of solutions to stochastic differential equations are also given. 



\bigskip
{\bf Plan of the paper and main results}

\smallskip

Let $\mu$ be a probability measure on $\mathbb R^n$ with density $e^{-V}$ where $V$ is $\mathcal C^2.$  
\smallskip

In Section~\ref{sec-dim-log-sobolev}, we propose a method based on the Borell-Brascamp-Lieb inequality to get dimensional logarithmic Sobolev inequalities in the spirit of the works~\cite{BL00,BL09} by S. Bobkov and M. Ledoux. The method is based on a general convexity inequality given in Theorem~\ref{thm-super-etoile}.   For instance, in Corollary~\ref{coro-lsbe} we shall prove the following :
If $\He{V}\geq R\,{\rm Id}_n$ with $R>0$, then
\begin{equation}\label{LSn-intro}
\entf{\mu}(f^2)\leq n(s -1 - \log s )+\frac{1}{2R}\int \Big|(1-s)\nabla V+2s\frac{\nabla f}{f}\Big|^2 f^2 \, d\mu
\end{equation}
for any $s>0$ and any function $f$ such that $\int f^2d\mu=1$.
This improves upon the classical logarithmic Sobolev inequality~\eqref{LS} under the Bakry-\'Emery condition, which is recovered for $s=1$. 
\medskip

In Section~\ref{sec-talagrand} (Theorem~\ref{thm-tala}) we propose a dimensional Talagrand inequality through optimal transportation in the spirit of Barthe-Kolesnikov~\cite{BK08} and D. Cordero-Erausquin~\cite{cordero} or the recent~\cite{cordero15} : If $\He{V}\geq R\,{\rm Id}_n$ with $R>0$ then
\begin{equation}\label{T-intro}
\frac{R}{2} W_2^2(\mu, \nu) \leq \nu (V) - \mu (V)  + n - n \exp \Big[\frac{1}{n} \Big(\nu(V) - \mu (V) - H(\nu \vert \mu) \Big) \Big]
\end{equation}
 for all $\nu\in P_2(\R^n)$. This bound implies the classical Talagrand inequality~\eqref{talaclass}.  
Let us observe that, using the terminology of the $\Gamma_2$-condition, the associated Markov generator $L = \Delta - \nabla V \cdot \nabla $ does not satisfy a $CD(R,n)$ curvature dimension condition, but only $CD(R, \infty)$. In particular the general dimensional log Sobolev or Talagrand inequalities, obtained on manifolds (see~\cite{bgl-book}) or on abstract measure spaces (as in~\cite{EKS13}) do not hold. In Section~\ref{sec-concentration} we show how the dimensional corrective term in our new Talagrand inequality enables to get sharp concentration inequalities. 

\medskip
Inspired by recent results on the equivalence between contraction and $CD(R,n)$ condition in abstract measure spaces (see \cite{ambrosio-gigli-savare, EKS13, BGGK}), in Section~\ref{sec-application-fp}  we consider applications to refined dimensional contraction properties under $CD(R,\infty)$ (see Proposition~\ref{propcontr} and Corollary~\ref{coro-cvOu}); we shall see how the dimension improves the asymptotic behaviour for the laws of solutions to stochastic differential equations  (in the spirit of \cite{BGG11,BGG13}).  Again the generator $L = \Delta - \nabla V \cdot \nabla $ does not satisfy a $CD(R,n)$ condition, but only $CD(R, \infty)$. The key point here is to take advantage of the contribution of the diffusion term, which includes a dimensional term. We shall also see how the dimension influences the short time smoothing effect, through very simple arguments (see Proposition~\ref{t=0}). 

\medskip

In section~\ref{sec-dim-BL-ine} we prove two kinds of dimensional Brascamp-Lieb inequalities, a first one by a $L^2$ argument, a second one by a linearization argument in the Borell-Brascamp-Lieb inequality. For instance, under the sole assumption~$\He{V}>0$, Theorem~\ref{thm-brascamp-lieb-2} states that
\begin{equation}\label{BL-intro}
\varf{\mu}(f) \leq \int \nabla f \cdot \He{V}^{-1} \, \nabla f \, d\mu - 
\int \frac{(f-\nabla f \cdot\He{V}^{-1} \, \nabla V)^2}{n+{\nabla V \cdot \He{V}^{-1} \, \nabla V}}d\mu
\end{equation}
for any smooth function $f$ such that $\int fd\mu=0.$ We shall discuss the optimality of our bounds and compare them with other very recent dimensional refinements of the Brascamp-Lieb inequality.

\medskip
In the Gaussian case where $\mu = \gamma$, then the logarithmic Sobolev~\eqref{LSn-intro} (by optimising over $s$) and Talagrand~\eqref{T-intro} inequalities are exactly~\eqref{eq-logsob-bakry-ledoux} and~\eqref{taladimgauss} respectively, while the Poincar\'e inequality~\eqref{BL-intro} improves upon~\eqref{Pdimgauss}.

\bigskip

{\it Notation}: whenever there is no ambiguity we shall respectively use $H, I,  W_2, \delta_{LSI}$ and $\delta_{Tal}$ for $H(\nu|\mu)$, $I(\nu|\mu), W_2(\nu,\mu), \delta_{LSI}(\nu|\mu)$ and $\delta_{Tal}(\nu|\mu)$. We shall sometimes let $\entf{dx} (f) = \int f \log f dx$ and $\mu(f)=\int fd\mu$ and use the same notation for an absolutely continuous measure with respect to Lebesgue measure, and its density.


\section{Logarithmic Sobolev inequalities}
\label{sec-dim-log-sobolev}

The Pr\'ekopa-Leindler inequality is a reverse form of the H\"older inequality.  Let $F$, $G$, $H$ be non-negative measurable functions on $\dR^n$ satisfying $\int Fdx=\int Gdx=1$, and let $s,t\geq0$ be fixed such that $t+s=1$. Under the hypothesis
\begin{equation}\label{eq-PL}
H(t x + s y) \geq F(x)^t G(y)^s
\end{equation}
for any $x,y \in \dR ^n$, the Pr\'ekopa-Leindler inequality ensures that 
$\int H dx \geq 1$, see~\cite[Chap.~19]{villani-book1} for instance.

The Borell-Brascamp-Lieb inequality  is a stronger and dimensional form of  the Pr\'ekopa-Leindler inequality. 
Assume again $\int Fdx=\int Gdx=1$ and in addition that $F$, $G$ and  $H$ are positive; then  the Borell-Brascamp-Lieb inequality asserts  that $\int H dx \geq 1$ as soon as 
\begin{equation}
\label{eq-borell-bl}
H(tx+sy)\geq \Big( tF(x)^{-1/n}+sG(y)^{-1/n} \Big)^{-n} 
\end{equation}
for any $x,y\in \R^n$, instead of the stronger~\eqref{eq-PL} (by convexity); see again~\cite{villani-book1}. 

\medskip

The Pr\'ekopa-Leindler inequality in particular implies many geometrical and functional inequalities as logarithmic Sobolev and  Brascamp-Lieb inequalities, as observed by S.~Bobkov and M.~Ledoux in~\cite{BL00,BL09} (see also~\cite{gentil08} for an application to the modified logarithmic Sobolev inequality). In the coming sections we shall see how the Borell-Brascamp-Lieb inequality implies dimensional form of these inequalities. Following S.~Bobkov and M.~Ledoux~\cite{BL00,BL09} our proofs are based on Taylor expansions when $s\rightarrow 0$ or $F\rightarrow 0$.

\subsection{A general  convexity inequality via the Borell-Brascamp-Lieb inequality}
\label{sec-convexe-ine}

Let us first state a general consequence of the Borell-Brascamp-Lieb inequality. It will lead to various dimensional logarithmic Sobolev inequalities. 

In the sequel  we let $\psi^*$ be the Legendre transform of a function $\psi$ on $\mathbb R^n$, defined for $y \in\R^n$ by 
$$
\psi^*(y)=\sup_{x\in\R^n}\{y\cdot x-\psi(x) \}\in(-\infty,+\infty].
$$ 
If $\psi$ is $C^1$ and strictly convex satisfying  
$$
\lim_{|x|\rightarrow + \infty}\frac{\psi(x)}{|x|}=+\infty,
$$
then (see~\cite[Sect.~2.1.3 and~2.4.3]{villani-otp} for instance) for all $x \in \mathbb R^n$, $\psi^*(x)\in\R$ and 
\begin{equation}\label{legendre}
\psi(x) = \nabla \psi (x) \cdot x - \psi^* (\nabla \psi (x))
\quad
\textrm{and}
\quad
\nabla \psi^* (\nabla \psi(x)) = x.
\end{equation}

\begin{thm}[Convexity inequality]
\label{thm-super-etoile}
Let $g,W$ be  $\mathcal C^1$ and positive functions on $\R^n$ satsifying the normalization condition  $\int g^{-n}dx=\int W^{-n}dx=1$. Assume moreover that there exists a constant $C>0$ such that for all $x\in\R^n$, 
\begin{equation}
\tag{H1}
W(x)\geq \frac{1}{C}\frac{|x|^2}{2},
\end{equation}
\begin{equation}
\tag{H2}
\frac{1}{C}(|x|^2+1)\leq g(x)\leq C(|x|^2+1)\quad{\rm and}\quad |\nabla g(x)|\leq C(|x|+1).
\end{equation}
Then 
\begin{equation}
\label{eq-wetoile}
\int \frac{W^*(\nabla g)}{g^{n+1}}dx\geq0.
\end{equation}
If $W$ is a $\mathcal C^1$ positive and strictly convex function which satisfies (H1) and $\int W^{-n}dx=1$, then~\eqref{eq-wetoile} is an equality for $g=W$. 
\end{thm}

 The same statement can be proved for a larger class of functions $g$ and $W$. We only state this result with these restrictive hypotheses for simplicity reasons, as this setting will be sufficient for our main  application.

The rigorous proof is postponed to the Appendix~\ref{sec-appendix-1}. The idea is to perform a Taylor expansion of the  Borell-Brascamp-Lieb inequality~\eqref{eq-borell-bl} when $s=1-t$ goes to 0.
Indeed, let $F=g^{-n}$ and $G=W^{-n}$ in~\eqref{eq-borell-bl}, hence satisfying $\int Fdx=\int Gdx=1$. Then the function $H_t$ defined~by
\begin{equation}\label{eq-defHt}
H_t(z)^{-1/n}=\inf_{h\in\R^n}\BRA{tg\left(z+\frac{s}{t}h\right)+sW(z-h) }
\end{equation}
for $z\in\R^n$ satisfies $\int H_tdx\geq1$.  The first-order Taylor expansion of $H_t$, when $s=1-t$ goes to 0, gives 
$$
H_t(z)=g(z)^{-n}-s \, n \, g(z)^{-n-1} \big( z\cdot\nabla g(z)-g(z) \big) + s \, n \, \frac{W^*(\nabla g (z))}{g^{n+1} (z)}+o(s).
$$
Since
$$
\int g^{-n-1} ( z\cdot\nabla g-g) \, dx=0 
$$
by integration by parts, the Taylor expansion of $\int H_tdx\geq 1$ implies the inequality \eqref{eq-wetoile}.

\medskip

Applications of Theorem~\ref{thm-super-etoile} are described in the coming two sections. They are based on the following observation.  Let $V$ be a given function and let $W=e^{\frac{V}{n}}$. Then, from the convexity of the exponential function, for any $a\in\R$ and $y\in\R^n$,
\begin{equation*}
W^*(y)\leq \frac{1}{n}e^{a}V^*(ne^{-a} y)+(a-1)e^a.
\end{equation*}
Combined with~Theorem~\ref{thm-super-etoile}, this gives the following corollary which is the main tool in our applications:
\begin{coro}
\label{cor-apres-theo}
Under the hypotheses of Theorem~\ref{thm-super-etoile}, let $V=n\log W$. Then for any function $a$, 
\begin{equation}
\label{eq-suite}
\int \frac{1}{g^{n+1}}\PAR{e^{a}V^*(ne^{-a} \nabla g)+n(a-1)e^a} \, dx\geq0. 
\end{equation}
\end{coro}

\subsection{Euclidean logarithmic Sobolev inequalities}
\label{sec-euclidean-logsob}

As a warm up, let us first see how to quickly recover the classical Euclidean logarithmic Sobolev inequality, using \eqref{eq-suite}.
Let $C :\dR^n\rightarrow\dR^+$ be a strictly convex function such that $\int e^{-C}dx<+\infty$, and let us apply~\eqref{eq-suite} with $V=C+\beta$ and $W=e^{V/n}$; here $\beta=\log\int e^{-C}dx$ so that $\int e^{-V} dx = 1$. Since $V$ is convex and $\int e^{-V}dx<+\infty$, it is classical that $V$ grows at least linearly at infinity, so that $W$ satisfies hypothesis~(H1) .

Then let $p>1$. Let also $f$ be a $\mathcal C^1$ positive function such that $\int f^p  dx = 1$ and $g=f^{-p/n}$  satisfies~(H2), and let $a=-\frac{p}{n}\log f+u$ where $u$ is a real constant. Then $V^* = C^* - \beta$ and~\eqref{eq-suite} can be written as
\begin{equation}
\label{eq-logsob-a}
\forall u\in\mathbb R, \quad\int f^p\log(f^p) \, dx\leq n(u-1)-\beta+\int C^*\PAR{-p e^{-u}\frac{ \nabla f}{f}}f^p \, dx.
\end{equation}


We can optimise over $u$ in $\mathbb R$ in the following case. 
Suppose that
there exists $q>1$ such that $C$ is $q$-homogeneous, that is, 
$C(\lambda x)=\lambda ^q C(x)$ for any $\lambda\geq0$ and $x$ in $\mathbb R^n$. Then $C^*$ is $p$-homogeneous with $1/p+1/q=1$, and in particular above $C^*\PAR{-p e^{-u}{ \nabla f}/{f}}= p^pe^{-pu}f^{-p} C^*(- \nabla f)$. Thus inequality~\eqref{eq-logsob-a}~gives 
\begin{equation}
\label{eq-logsob-b}
\int f^p\log(f^p) \, dx\leq n(u-1)-\beta+e^{-pu}p^p\int C^*\PAR{ - \nabla f} \, dx
\end{equation}
for any function $f$ such that $\int f^p dx = 1$ and $f^{-n/p}$ satisfies~(H2). Now, let $f$ be a $\mathcal C^1$ non negative and compactly supported function and for $\varepsilon>0$ let $f_\varepsilon(x)=C_{\varepsilon}(\varepsilon(|x|^2+1)^{-n/p}+f)$, where $C_\varepsilon$ is such that $\int (f_\varepsilon)^p dx=1$. The function $f_\varepsilon^{-n/p}$ satisfies~(H2) for any $\varepsilon$. Taking the limit when $\varepsilon$ goes to~$0$, inequality~\eqref{eq-logsob-b} then holds for any $\mathcal C^1$ non negative and compactly supported function $f$ such that $\int f^pdx=1$.

For the optimal $ u = p^{-1} \log \, ( \, p^{p+1} \int C^*(-\nabla f) dx / n ),$ the bound~\eqref{eq-logsob-b} leads to 
$$
\int f^p \, \log (f^p) \, dx \leq \frac{n}{p}\log\PAR{\frac{p^{p+1}}{ne^{p-1}}\frac{\int C^*(-\nabla f)dx}{(\int e^{-C}dx)^{p/n}}}
$$
for any $\mathcal C^1$ non negative and compactly supported function $f$ such that $\int f^pdx=1$. Of course, the inequality can be extended to a larger class of functions $f$. 
Hence, we recover the optimal $L^p$-Euclidean log Sobolev inequality proved in~\cite{del-pino-dolbeault03,gentil03}  and in particular, setting $C(x) = \vert x \vert^2/2$ and $p=q=2$, the classical inequality
$$
\int f^2 \, \log (f^2) \, dx \leq \frac{n}{2}\log\PAR{\frac{2}{n \pi e} \int \vert \nabla f \vert^2 dx}.
$$

\subsection{Dimensional logarithmic Sobolev inequalities}
\label{sec-logsob-2}

In this section we consider a probability measure $\mu$ with density $e^{-V}$ and the function $W=e^{V/n}$, and a positive function $f$ such that $\int f^2 \, d\mu =1$. We assume again that $V$ is convex ; then $W=e^{V/n}$ satisfies hypothesis~(H1) since  $\int e^{-V}dx=1$.

Corollary~\ref{cor-apres-theo} applied with $g=e^{V/n}\,f^{-2/n}$ (assuming that $g$ satisfies hypothesis~(H2)) and $a=\frac{V}{n}-\frac{2}{n}\log f+u$ with $u \in \mathbb R$ gives
$$
\int \PAR{V^*\Big(e^{-u} \nabla V-2e^{-u}\frac{\nabla f}{f}\Big)+V-\log(f^2)+n(u-1)} f^2  e^{-V}dx\geq 0.
$$
\begin{coro}
\label{cor-un-de-plus}
Let $d\mu(x)=e^{-V(x)}dx$ be a probability measure with $V$ a convex function and let $f$ be a $\mathcal C^1$ positive function such that $\int f^2 \, d\mu =1$ and such that $g=e^{V/n}\,f^{-2/n}$  satisfies hypothesis~(H2). 
Then for any $s>0$
\begin{equation}
\label{etoileetoile}
\entf{\mu}(f^2)\leq \int\SBRA{ V^*\Big(s\nabla V-2s\frac{\nabla f}{f}\Big)+ V}f^2d\mu - n(1 + \log s ).
\end{equation}
\end{coro}

For $s=1$, inequality~\eqref{etoileetoile} simplifies as 
$$
\entf{\mu}(f^2)\leq \int\SBRA{V^*\Big(\nabla V-2\frac{\nabla f}{f}\Big)+V-n}f^2d\mu, \qquad \int f^2 d\mu = 1.
$$
In particular, for $V=\frac{|x|^2}{2}+ \frac{n}{2} \log(2\pi)$, then $\mu$ is the standard Gaussian measure $\gamma$ and we recover the Gaussian logarithmic  Sobolev inequality of L.~Gross, 
$$
\entf{\gamma}(f^2)\leq {2}\int |\nabla f|^2d\gamma, \qquad \int f^2 \, d\gamma =1.
$$
More generally, let $V$ be a strictly convex function on $\R^n$. Then inequality~\eqref{etoileetoile} with $s=1$, by~\eqref{legendre} and integration by parts, leads to the modified logarithmic Sobolev inequality
$$
\entf{\mu}(f^2)\leq \int\SBRA{ V^*\Big(\nabla V-2\frac{\nabla f}{f}\Big)+2x\cdot\frac{\nabla f}{f}-V^*(\nabla V)}f^2d\mu, \qquad \int f^2 \, d\mu =1
$$
proved  by the second author in~\cite{gentil08}.

\bigskip

Assuming {\it uniform convexity} on $V$ we now optimise over the parameter $s>0$ in Corollary~\ref{cor-un-de-plus},  to obtain dimensional logarithmic Sobolev inequalities. 
 Suppose that $V$ is $\mathcal C^2$ with $\He{V} \geq R\,{\rm Id}_n$ for $R>0.$ Then, for their inverse matrices, $\He{V^*} \leq R^{-1} \,{\rm Id}_n$ on $\mathbb R^n$. Hence, for any $z$ and by the Taylor expansion at point $\nabla V(x)$, 
\begin{multline*}
V^*(z) + V (x) \leq V^*(\nabla V(x)) + \nabla V^*(\nabla V(x)) \cdot (z - \nabla V(x)) + \frac{1}{2R} \vert z - \nabla V(x) \vert^2 + V (x)\\
= x \cdot z + \frac{1}{2R} \vert z - \nabla V (x) \vert^2.
\end{multline*}
Here we use the relations~\eqref{legendre}. 
For $z = s  \nabla V-2s\frac{\nabla f}{f}$ at point $x,$ and by~\eqref{etoileetoile}, this leads to
$$
\entf{\mu}(f^2)\leq -n(1 + \log s)+s \int x \cdot \Big(\nabla V-2\frac{\nabla f}{f}\Big) \, f^2d\mu + \frac{1}{2R}\int \Big|\Big(s\nabla V-2s\frac{\nabla f}{f}\Big) - \nabla V\Big|^2 f^2d\mu.
$$
By integration by parts and extending to compactly supported functions, as for~\eqref{eq-logsob-b}, we finally obtain:

\begin{coro}[Dimensional LSI under $\Gamma_2$-condition]
\label{coro-lsbe}
Let $\mu$ be a probabili\-ty mea\-sure with density $e^{-V}$ where $V$ is $\mathcal C^2$  with $\He{V}\geq R\,{\rm Id}_n$ for $R>0$. Then
\begin{equation}
\label{eq-log-be}
\entf{\mu}(f^2)\leq n(s -1 - \log s )+\frac{1}{2R}\int \Big|(1-s)\nabla V+2s\frac{\nabla f}{f}\Big|^2 f^2 \, d\mu
\end{equation}
for any $s>0$ and any $\mathcal C^1$, non-negative and compactly supported function $f$ such that $\int f^2d\mu=1$.
\end{coro}

The bound can of course be extended to other classes of functions $f$.

When $s=1$, we recover the classical logarithmic Sobolev inequality ~\eqref{LS} under the Bakry-\'Emery condition. 

%

Let us observe that the right-hand side in~\eqref{eq-log-be} can be expanded as $-n \log s$ plus a second order polynomial in $s$. Hence it admits a unique minimiser $s>0$, which solves a second order polynomial. The obtained expression is not appealing and we prefer to omit it.  
In the Gaussian case where $\mu = \gamma$, then the optimisation over $s$ gets even simpler and leads again to the dimensional Gaussian log Sobolev inequality~\eqref{eq-logsob-bakry-ledoux}. 

Moreover, for a general $V$ and as in~\eqref{defLSI} or~\eqref{deftala} below for the Talagrand inequality, the bound~\eqref{eq-log-be} can be written as a (not either appealing) deficit in the log Sobolev inequality. 

\medskip

We will see in Section \ref{regu} that \eqref{eq-log-be} leads to new and sharp short time smoothing on the entropy of solutions to an associated Fokker-Planck equation.


\section{Talagrand inequalities}
\label{sec-talagrand}

 The main result of this section is 

\begin{thm}\label{thm-tala} {\bf (Dimensional Talagrand inequality)}
Let $\mu$ be a probability measure in $P_2(\R^n)$ with density $e^{-V}$ where $V$ is a  $\mathcal C^2$ function satisfying $\He{V} \geq R\,{\rm Id}_n$ with $R>0$. Then for all $\nu\in P_2(\R^n)$
\begin{equation}\label{tala}
\frac{R}{2} W_2^2(\mu, \nu) \leq \nu (V) - \mu (V)  + n - n \exp \Big[\frac{1}{n} \Big(\nu(V) - \mu (V) - H(\nu \vert \mu) \Big) \Big].
\end{equation}
\end{thm}

\medskip

In other words, if $\He{V} \geq R\,{\rm Id}_n$, then $\nu(V) - \mu(V) - \frac{R}{2} W_2^2(\nu, \mu) >-n$ and
\begin{equation}\label{deftala}
\delta_{Tal} (\nu \vert \mu) \geq  \max 
\Big\{ 
\delta_n \Big( H(\nu|\mu)  +\mu (V) -\nu (V))\Big), 
\Lambda_n \Big(\nu(V) - \mu(V) - \frac{R}{2} W_2^2(\nu, \mu) 
\Big) 
\Big\}.
\end{equation}
Here $\delta_n$ and $\Lambda_n$ are the positive functions respectively defined by $\delta_n(x) = n [ e^{-x/n} - 1 + x/n ], x \in \mathbb R$ and $\Lambda_n(x) = x - n \log (1+x/n), x > -n.$

\smallskip

The function $\delta_1(x)=e^{-x}-1+x$ is positive and convex. It is moreover decreasing on $\mathbb R^-$ and increasing on $\mathbb R^+$. By a direct computation, $\delta_1(x)$ is bounded from below by $x^2/2$ if $x \leq 0$, $x^2/e$ if $0 \leq x \leq 1$ and $x/e$ if $x >1$; hence always by $\frac{1}{e} \min(|x|,x^2)$. Then for any $x\in\R$, $\delta_n(x) \geq \frac{1}{e} \min(|x|,\frac{x^2}{n}).$

\medskip

Since $e^u \geq 1+u,$ the bound \eqref{tala} implies the classical Talagrand inequality~\eqref{talaclass} under the condition $\He{V} \geq R\,{\rm Id}_n$.  When $\mu$ is the standard Gaussian measure $\gamma$ on~$\rr^n$, then $R=1$ and we recover the dimensional Talagrand inequality~\eqref{taladimgauss}.

\medskip

Under a moment condition Theorem~\ref{thm-tala} simplifies as follows:

\begin{coro}
Following the same assumptions as in Theorem~\ref{thm-tala},  for all~$\nu$ in $P_2(\R^n)$ with  $\nu (V) \le \mu(V)$, 
\begin{equation}\label{tala-c}
\delta_{Tal} (\nu \vert \mu) \ge \delta_n(H(\nu|\mu)) \ge \frac{1}{e} \min\left(H(\nu|\mu), \frac{H(\nu|\mu)^2}n\right).
\end{equation}
\end{coro}

Theorem \ref{thm-tala} will be deduced from the following dimensional $HWI$-type inequality, applied with $f=1$ and $\nu = g \mu$. The $HWI$ inequality bounds from above the entropy by the Wasserstein distance and the Fisher information (defined in~\eqref{eq-def-fisher}), in the form
\begin{equation}\label{hwi0}
H(\nu|\mu) \leq W_2(\nu, \mu) \, \sqrt{I(\nu \vert \mu)} - \frac{R}{2} W_2^2(\mu, \nu)
\end{equation}
for all $\nu$. It has been introduced in~\cite{ov00} and proved in~\cite{ov00} and~\cite{cordero} under the Bakry-\'Emery condition $\He{V} \geq R\,{\rm Id}_n, R \in \mathbb R.$
\begin{thm}[Dimensional  $HWI$ inequality]
\label{thm-hwi}
Let $\mu$ be a probability measure on $\R^n$ with density $e^{-V}$ where $V$ is a  $\mathcal C^2$ function satisfying $\He{V} \geq R\,{\rm Id}_n$ with $R \in \mathbb R$.  Let also $f,g$ be smooth functions such that $f\mu$ and $g\mu$ belong to $P_2(\R^n)$. Then
\begin{multline*}
 n \exp \Big[\frac{1}{n} \Big(H(f\mu|\mu) - H(g \mu \vert \mu) + \mu (gV)- \mu (fV)  \Big) \Big] - n\\
 \leq
  \mu (gV)- \mu (fV) + W_2(f\mu, g\mu) \sqrt{I(f\mu \vert \mu)} - \frac{R}{2} W_2^2(f\mu, g\mu).
\end{multline*}
\end{thm}

\smallskip

For $g= 1$ and $\nu = f \mu$, this bound can be written as the dimensional $HWI$ inequality
\begin{equation}
\label{hwi}
 n \exp \Big[\frac{1}{n} \Big(H(\nu \vert \mu) + \mu (V)- \nu (V) \Big) \Big] - n \leq
 \mu (V) - \nu (V)+ W_2(\mu, \nu) \sqrt{I(\nu \vert \mu)} - \frac{R}{2} W_2^2(\mu, \nu).
\end{equation}
As in \eqref{deftala} for the Talagrand inequality, this can equivalently be written as a deficit in the $HWI$ inequality. It is classical that the $HWI$ 
inequality~\eqref{hwi0} implies the logarithmic Sobolev inequality~\eqref{LS} (see~\cite{ov00} for instance). 
Likewise, from~\eqref{hwi}, one can obtain a dimension dependent logarithmic Sobolev inequality. We refer to Section~\ref{rem-LSidim} for further details.

 The proof of Theorem~\ref{thm-hwi} will be given in Section~\ref{thm-23}.

\subsection{An application to concentration}
\label{sec-concentration}
Let us quickly revisit K. Marton's argument for concentration via Talagrand's inequality (as in \cite[Chap.~22]{villani-book1} for instance) and see how the refined inequality~\eqref{tala} in Theorem \ref{thm-tala} gives sharpened information for large deviations.

\medskip

Let $d\mu = e^{-V} dx$ satisfy inequality~\eqref{tala}. Let also $A\subset \rr^n$, $r>0$ and $A_r=\{x;\,\forall y \in A,\, \vert y-x \vert >r\}$. Let finally $\mu_A=\frac{1_A}{\mu(A)}\mu$ and $\mu_{A_r}=\frac{1_{A_r}}{\mu(A_r)}\mu$ be the restrictions of $\mu$ to $A$ and $A_r$. Then, as $W_2$ is a distance,
$$
r\le W_2(\mu_A,\mu_{A_r})\le W_2(\mu_A,\mu)+W_2(\mu_{A_r},\mu).
$$
First of all 
$$
W_2(\mu_A,\mu) \leq \sqrt{2 R^{-1} H(\mu_A \vert \mu)} = \sqrt{2R^{-1}\log(1/\mu(A))} := c_A
$$ 
by \eqref{tala}, or its weaker form~\eqref{talaclass}. Let now $c_V=\int Vd\mu, x_r=H(\mu_{A_r}|\mu)=\log(1/\mu(A_r))$ and $V_r= \int V d\mu_{A_r}$. By \eqref{tala} again we get, for $r>c_A$, 
$$
(r-c_A)^2\le W_2^2(\mu_{A_r},\mu)\le \frac{2}{R}\Big( V_r-c_V+n-n\exp \Big[ -\frac{1}{n} (x_r+c_V-V_r)\Big]\Big).
$$
Since  $x_r=\log(1/\mu(A_r))$ we obtain :
\begin{coro}[Concentration inequality]
\label{coro-conc}
Following the same assumptions as in Theorem~\ref{thm-tala}, let $A\subset\rr^n$, $r>0$ and $A_r=\{x;\,\forall y \in A,\, \vert y-x \vert >r\}, c_A = \sqrt{2R^{-1}\log(1/\mu(A))}, c_V=\int Vd\mu, V_r=\int  V d\mu_{A_r}$. Then for $r>c_A$
$$
\mu(A_r) \leq e^{c_V- V_r} \Big[ 1 + \frac{1}{n} \big(V_r - c_V - \frac{R}{2} (r-c_A)^2 \big) \Big]^n.
$$
\end{coro}
Since $(1+u/n)^n \leq e^u$, the bound in Corollary \ref{coro-conc} implies the classical Gaussian concentration
$$
\mu(A_r) \leq e^{-\frac{R}{2} (r-c_A)^2}, \qquad r > c_A
$$
of the Talagrand inequality~\eqref{talaclass}, see again~\cite[Chap.~22]{villani-book1} for instance.


The bound in Corollary \ref{coro-conc} captures the behaviour of concentration of the measure $\mu$ in a more accurate way: let for instance $V(x) = \vert x \vert^2/2 + \vert x \vert^p +Z_p$ with $p > 2$ and a normalizing factor $Z_p$, and $A$ be the Euclidean unit ball in $\mathbb R^n$. Then $\He{V} \geq {\rm Id}_n$, so by Corollary \ref{coro-conc}  with $R=1$ there exists a constant $C = C(p,n)$ such that for all $r  > C$
$$
\mu(\vert x \vert > r+1) = \mu(A_{r}) \leq \exp \Big[ c_V - V_r + n \, \log(1+V_r/n) \Big].
$$ 
 But $V_r \geq r^p+Z_p$,  so for all $\varepsilon <1$ there exists another constant $C$ depending also on $\varepsilon$ such that~for~all~$r>C$
$$
\mu(\vert x \vert > r) \leq e^{-(1- \varepsilon) r^p}. 
$$ 
This concentration inequality in this precise example can also be obtained by using a $L^p$-Talagrand inequality or a $L^p$-log Sobolev inequality; however we have found it interesting to get it by means of the dimension dependence of the classical Talagrand inequality, moreover in a shorter and more straightforward manner.

\subsection{Tensorisation and comparison with earlier results}\label{tenso}

 In $\mathbb R^n,$ let $W_1$ be the Wasserstein distance between probability measures, for the cost $\vert y-x \vert, x, y \in \mathbb R^n.$

\smallskip

Deficit in the Gaussian Talagrand inequality (for $\mu = \gamma$) and for centered measures $\nu$ has been investigated in {\it one dimension} in~\cite{BK08} and~\cite{bgrs14}, in the form
$$
\delta_{Tal}(\nu \vert \gamma) \geq c \inf_{\pi} \int_{\mathbb  R \times \mathbb R} \Lambda(\vert y-x \vert) d\pi(x,y) \geq c \min \big\{W_1(\nu, \gamma)^2, W_1(\nu, \gamma) \big\}.
$$
Here the $c$'s are diverse numerical constants and the infimum runs over couplings $\pi$ of $\gamma$ and $\nu$.

This second lower bound has been extended in~\cite[Th. 5]{FIL14} to {\it any dimension $n$}, as
\begin{equation}\label{defW11}
\delta_{Tal} (\nu \vert \gamma) \ge c\min\left(\frac{W_{1,1}(\nu,\gamma)^2}n,\frac{W_{1,1}(\nu,\gamma)}{\sqrt{n}}\right)
\end{equation}
as soon as $\nu$ has mean $0$; here $c$ is a numerical constant independent of the dimension $n,$ and on $\mathbb R^n \times \mathbb R^n$
$$
W_{1,1}(\mu,\nu)=\inf_\pi \int_{\mathbb  R^n \times \mathbb R^n} \sum_{i=1}^n|y_i-x_i| \; d\pi(x,y).
$$

Still under a centering condition, the bound~\eqref{defW11} has been improved in~\cite[Prop. 3]{cordero15} by replacing the quantity $W_{1,1} / \sqrt{n}$ by the larger $W_1$ Wasserstein distance on $\mathbb R^n$, and extended to reference measures $\mu$ with density $e^{-V}$ where $\He{V} \geq R\,{\rm Id}_n.$

In comparison, our bound~\eqref{deftala} has the following two advantages : it holds without any centering condition on $\nu$, and gives a lower bound on the deficit in terms of the relative entropy $H$ : this is a strong way of measuring the gap between measures, by the Pinsker inequality for instance (see~\cite[Chap. 22]{villani-book1}), and the relative entropy can be much larger than the weak distance $W_2$.

\smallskip

As considered in~\cite{cordero15} and~\cite{FIL14}, a natural example is the product measure case when $\mu^N=\otimes^N\mu$ and $\nu^N =\otimes^N\nu$ on $\mathbb R^{nN}$ for $N\in\N^*$. Then $\delta_{Tal}(\nu^N|\mu^N)=N \, \delta_{Tal}(\nu|\mu)$ by tensorisation properties of both $H$ and $W_2^2$. However, the above bound \eqref{defW11} in~\cite{cordero15} (so with $W_1$ instead of $W_{1,1} / \sqrt{n}$) gives a lower bound on $\delta_{Tal}(\nu^N|\mu^N)$ equal to a constant $c$ times
$$
\min\left(W_1(\nu^N \! ,\mu^N)^2,W_1(\nu^N \!,\mu^N)\right)
\leq
\min\left(W_2(\nu^N \!,\mu^N)^2,W_2(\nu^N \!,\mu^N)\right)
=
\min\left(NW_2(\nu,\mu)^2,\sqrt{N}W_2(\nu,\mu)\right)
$$
since $W_1 \leq W_2$. Hence this lower bound has the good order in $N$ at most only for small perturbations $\nu$ of the reference measure $\mu$.

In contrast, our bound always has the correct order in $N$. Indeed, if  $V^{(N)} = \oplus^N V$ so that $d\mu^N=e^{-V^{(N)}}dx$ on $\mathbb R^{nN},$ then
$$
H(\nu^N|\mu^N)+\mu^N (V^{(N)}) - \nu^N (V^{(N)})=N\left(H(\nu|\mu)  +\mu(V) -  \nu (V)\right);
$$
hence Theorem \ref{thm-tala} leads to 
$$
\delta_{Tal}(\nu^N|\mu^N)\ge N\, \delta_n\left(H(\nu|\mu)+ \mu (V) - \nu(V)\right),
$$
which has the correct order in $N$. 

 \subsection{Useful facts on optimal transport}\label{rappels}

In the proof of Theorem \ref{thm-hwi} and in proofs below we shall need the following notation and facts. 

\smallskip

If $\mu$ is a probability measure on $\mathbb R^n$ and $T : \mathbb R^n \to \mathbb R^n$ a Borel function, we let $T \# \mu$ be the image measure of $\mu$ by $T$, defined by $T \# \mu (h) = \mu(h \circ T)$ for all bounded continuous functions $h : \mathbb R^n \to \mathbb R.$

Let now $\mu_0$ and $\mu_1$ in $P_2(\R^n)$ be absolutely continuous with respect to Lebesgue measure. Then there exists a convex function $\varphi$ on $\mathbb R^n$ such that $\mu_1 = \nabla \varphi\# \mu_0$ (see~\cite[Th.~2.12]{villani-otp} or \cite[Th.~10.41]{villani-book1} for instance). The map $\nabla \varphi$ is called the Brenier map. Moreover
$$
\int \vert \nabla \varphi(x) - x \vert^2 d\mu_0(x) = W_2^2(\mu_0, \mu_1).
$$

Now, by the Alexandrov Theorem (see~\cite{mccann97} or~\cite[Th.~14.1]{villani-book1} for instance), a convex function $\psi$ is almost everywhere twice differentiable: for almost every $x \in \mathbb R^n$ there exists a non negative symmetric matrix $A$ such that
$$
\psi(x+h) =\psi(x) + \nabla \psi(x) \cdot h + \frac{1}{2} A h \cdot h + o(\vert h \vert^2)
$$
as $h$ tends to $0$ in $\mathbb R^n.$ The matrix $A$ is denoted $\He{\psi} (x)$ and called the Hessian of $\psi$ in the sense of Alexandrov. The trace of $A$ will be denoted $\Delta \psi(x)$ : it coincides with the density of the absolutely continuous part of the distributional Laplacian of $\psi$, the singular part being a non negative measure. 

In fact, in the above notation and by~\cite[Th.~4.4]{mccann97} or~\cite[Th.~6.2.12]{ambrosio-gigli-savare},  $\He{\varphi} (x)$ is a positive matrix for $\mu_0$-almost every $x$. Moreover, by~\cite{mccann97} (see also~\cite[Lem.~5.5.3]{ambrosio-gigli-savare}), the Brenier map solves the Monge-Amp\`ere equation
\begin{equation}\label{MA}
\mu_0(x) = \mu_1(\nabla \varphi(x)) \det (\He{\varphi}(x))
\end{equation}
at $\mu_0$-almost every $x$ in $\R^n$. Here $\mu_0$ and $\mu_1$ are the densities of the measures.

\medskip

Let now $\varphi^*$ be the Legendre transform of $\varphi$. Then $\mu_0 = \nabla \varphi^* \# \mu_1$ by~\cite[Th.~2.12]{villani-otp} for instance. Moreover $\nabla \varphi^*(\nabla \varphi(x)) = x$ and $\nabla \varphi(\nabla \varphi^*(y)) = y$ for $\mu_0$-almost every $x$ and $\mu_1$-almost every $y$. 

Furthermore, by~\cite[Th.~A.1]{mccann97}, if $\He{\varphi} (x)$ is invertible at $x$ then $\varphi^*$ is twice differentiable at $\nabla \varphi (x),$ with
$\He{\varphi^*} (\nabla \varphi (x)) = \big[ \He{\varphi} (x) \big]^{-1}.$ By the remark above, this is the case for $\mu_0$-almost every $x.$

\medskip

Finally, the curve $(\mu_s)_{s \in [0, 1]}$ defined by $\mu_s = ((1-s) Id + s \nabla \varphi) \# \mu_0$ is a geodesic path in $P_2(\mathbb R^n)$ between $\mu_0$ and $\mu_1$, in the sense that 
$$
W_2(\mu_s, \mu_t) = \vert t-s \vert \; W_2(\mu_0, \mu_1)
$$
for all $0 \leq s, t \leq 1.$ It holds that $\mu_s$ is also absolutely continuous with respect to Lebesgue measure, see~\cite[Prop.~1.3]{mccann97} or~\cite[Th.~5.9]{villani-otp} for instance.

\subsection{Proof of Theorem \ref{thm-hwi}}
\label{thm-23}

Theorem~\ref{thm-hwi} is a consequence of the relation  
\begin{equation}\label{entropies}
H(h\mu|\mu) - \mu (hV) = \entf{dx} (h e^{-V})
\end{equation}
written with $h =f, g$ and of the following lemma. 
\begin{lemma}
\label{lemma-fg}
Following the same assumptions as in Theorem~\ref{thm-hwi}, let $f,g$ be two smooth functions such that $f\mu$ and $g\mu$ belong to $P_2(\R^n)$. Let $\varphi$ be a convex function on $\rr^n$ such that 
$\nabla \varphi \# (f \mu) = g \mu$. Then
\begin{multline*}
 \int V  \, g \, d\mu - \int V \, f \, d\mu - \int (\nabla \varphi -x) \cdot \nabla f \, d\mu 
 \geq
   n \exp \Big[\frac{1}{n} \Big(\entf{dx} (f e^{-V}) - \entf{dx} (g e^{-V}) \Big) \Big] - n
 \\  
+    \int  \int_0^1 (\nabla \varphi(x) - x) \cdot \He{V} (x + t (\nabla \varphi(x) -x))  (\nabla \varphi(x) - x) (1-t) dt \, f(x) \, d\mu(x).
\end{multline*}
\end{lemma}

Indeed, if $\He{V} \geq R\,{\rm Id}_n$, then the last term above is greater than $\frac{R}{2} \int \vert \nabla \varphi -x \vert^2 f \, d\mu = \frac{R}{2} W_2^2(f\mu, g\mu)$. Moreover, on the left-hand side,
$$
-  \int (\nabla \varphi -x) \cdot \nabla f \, d\mu \leq \Big[ \int \vert \nabla \varphi -x \vert^2 \, f \, d\mu \Big]^{1/2} \Big[ \int \frac{\vert \nabla f \vert^2}{f} d \mu \Big]^{1/2} = W_2(f\mu, g\mu) \sqrt{I (f \mu \vert \mu)}
$$
by the Cauchy-Schwarz inequality.  This implies Theorem~\ref{thm-hwi}.
\bigskip

{\bf Proof of Lemma~\ref{lemma-fg}.}

 By the Taylor formula,
$$
V (\nabla \varphi(x)) - V(x) = \nabla V(x) \cdot (\nabla \varphi(x) - x) + \int_0^1 (\nabla \varphi(x) - x) \cdot \He{V} (x + t (\nabla \varphi(x) -x))  (\nabla \varphi(x) - x) (1-t) dt
$$ 
for almost every $x$ in $\mathbb R^n$. We now  integrate with respect to $f \, \mu$ and use the comparison between Alexandrov and distributional Laplacians to deduce that
$$
\int \nabla V(x) \cdot (\nabla \varphi(x) - x) \, f(x) \, d\mu(x) \geq \int \big[ (\Delta \varphi - n) f + (\nabla \varphi -x) \cdot \nabla f \big] \, d\mu = \int \Delta \varphi \, f \, d\mu -n + \int (\nabla \varphi -x) \cdot \nabla f \, d\mu,
$$
as in~\cite{cordero} or~\cite[Th.~9.17]{villani-otp} for instance.
This leads~to
\begin{multline}\label{Vconv}
\int V \, g \, d\mu - \int V \, f \, d\mu - \int (\nabla \varphi -x) \cdot \nabla f \, d\mu \geq
 \int \Delta \varphi \,f\, d\mu -n \\
 + \int  \int_0^1 (\nabla \varphi(x) - x) \cdot \He{V} (x + t (\nabla \varphi(x) -x))  (\nabla \varphi(x) - x) (1-t) dt \, f(x) \, d\mu(x).
\end{multline}
Then Lemma \ref{lemma-fg} is a consequence of the following Lemma. 

\begin{lemma}\label{lemma-laplace}
Let $\mu_0,\mu_1\in P_2(\R^n)$ absolutely continuous with respect to Lebesgue measure, with respective densities also denoted $\mu_0$ and $\mu_1$.
Let $\varphi$ be a convex function on $\rr^n$ such that $\nabla \varphi \# \mu_0 = \mu_1$. Then 
\begin{equation}\label{bornelaplace}
\int \Delta \varphi \, d\mu_0 \geq n \, \exp \left[ \frac{\entf{dx} (\mu_0) - \entf{dx} (\mu_1)}{n} \right].
\end{equation}
\end{lemma}

\begin{eproof} 
Taking logarithms in the Monge-Amp\`ere equation~\eqref{MA} and integrating with respect to $\mu_0$ lead~to
\begin{equation}\label{entdx}
\entf{dx} (\mu_0) = \entf{dx} (\mu_1) + \int \log \det (\He{\varphi}) \, d\mu_0.
\end{equation}

Now, if for each $x$ the symmetric matrix $\He{\varphi}$ has eigenvalues $\varphi_i$, then by the Jensen inequality
$$
\int\! \log \det (\He{\varphi}) \, d\mu_0\! =\! n \, \frac{1}{n} \sum_{i} \int\! \log (\varphi_i) \, d\mu_0 \leq n \log \Big(\int\! \frac{1}{n} \sum_i \varphi_i \, d\mu_0 \Big) \!=\! n \log \Big(\frac{1}{n} \int\! \Delta \varphi  \, d\mu_0 \Big).
$$
This concludes the proof.
\end{eproof}

\begin{remark} 
 In the Gaussian case, we have already observed that translations of the Gaussian measure are extremals of the Talagrand inequality. As observed in~\cite{cordero}, or as can be observed from the proof above, there are no other extremals. Indeed the Hessian of the map $\varphi$ has to be constant and equal to the identity matrix for all inequalities to be equalities. 
 
 In fact, if $\He{V}\ge R\,{\rm Id}_n$, then equality in the Talagrand inequality implies that the potential is necessarily Gaussian and that extremals are translations of the Gaussian measure. 
\end{remark}

\subsection{Logarithmic Sobolev inequalities by transport}
\label{rem-LSidim}

As observed in~\cite{ov00}, the $HWI$  inequality~\eqref{hwi0} classically implies the logarithmic  Sobolev inequality~\eqref{LS} by bounding from above the second order polynomial in $W_2$ in $HWI$ by its maximum. Likewise, the dimensional $HWI$  inequality~\eqref{hwi} is another path towards dimensional logarithmic Sobolev inequalities. Here we obtain : 

Let $\mu$ have density $e^{-V}$ where $V$ is $\mathcal C^2$ and satisfies $\He{V} \geq R\,{\rm Id}_n$ with $R>0$. Then
$$
H(\nu \vert \mu) \leq
\nu(V) - \mu (V)+ n \log \Big( 1 + \frac{1}{n} \Big(\frac{I(\nu \vert \mu)}{2R} + \mu (V)- \nu (V) \Big) \Big)
$$
for all $\nu$.  Equivalently, in terms of deficit, 
\begin{equation}\label{defLSI}
\delta_{LSI} (\nu \vert \mu) \geq R \max \Big\{  \delta_n \Big( \nu (V)- \mu (V)- H(\nu \vert \mu) \Big), \Lambda_n \Big(\frac{I(\nu \vert \mu)}{2R} - \nu (V) + \mu(V)) \Big) \Big\}.
\end{equation}

In the Gaussian case, then $R=1$ and we obtain a bound which is slightly worse than~\eqref{eq-logsob-bakry-ledoux}, where a $\log(1+2u)$ term is replaced by the larger $2 \log(1+u)$.

At this point, let us observe that still in the Gaussian case a dimensional $HWI$  has been derived in \cite[Th.~1.1]{bgrs14}. It is also observed by the authors that the $HWI$ inequality in~\cite{bgrs14} does not seem to imply~\eqref{eq-logsob-bakry-ledoux}. 
We could not compare the $HWI$ in~\cite{bgrs14} to our bound~\eqref{hwi} in full generality. However, if $\nu( \vert x \vert^2)  = n = \gamma( \vert x \vert^2)$ then they can respectively be written  as
$$
2h \leq x-y + \log(1+x) \qquad \textrm{and} \qquad h \leq \log (1+x-y/2)
$$
for $x = W_2 \sqrt{I}/n, y = W_2^2/n$ and $h = H/n;$ hence our bound is at least significantly more precise in the common range $I \gg W_2 \sim 1$: indeed then $x \gg y \sim 1$ in this range, so that comparing the two right-hand sides amounts to $x \gg \log(1+x).$


\medskip

As remarked in \cite{bgrs14,FIL14} it is also possible to get refined logarithmic Sobolev inequalities by combining the $HWI$ and Talagrand inequalities. Here, if $\He{V} \ge R\,{\rm Id}_n$ with $R>0$, then~\eqref{hwi} can be written as
\begin{equation}\label{def-HW2}
H +  \delta_n(-h) \leq W_2 \sqrt{I} - \frac{R}{2} W_2^2
\end{equation}
where $h = H + \mu (V)- \nu (V)$. Moreover $H = \frac{R}{2} W_2^2+ \delta_{Tal}$, so
$$
\frac{\delta_{Tal} + \delta_n(-h)}{W_2} \leq \sqrt{I} - R W_2.
$$
Then, by~\eqref{def-HW2} again and Theorem~\ref{thm-tala},
\begin{eqnarray*}
\delta_{LSI} 
= \frac 12 I - R\,H
\ge
R \, \delta_n(-h) + \frac12\left(\sqrt{I}-  R\,W_2\right)^2 
&\ge& 
R \delta_n(-h) + \frac12 \frac{(\delta_{Tal} + \delta_n(-h))^2}{W_2^2}
\\
&\ge&
 R \delta_n(-h) + \frac12 \frac{(\delta_n(h) + \delta_n(-h))^2}{W_2^2}.
\end{eqnarray*}

In particular this improves upon the first lower bound in~\eqref{defLSI}.
Let us recall that the function $\delta_n$ is defined above, after Theorem~\ref{thm-tala}.

\medskip

Refined Gaussian logarithmic  Sobolev inequalities have been considered for certain classes of test measures~$\nu$ : measures $\nu$ satisfying lower and upper curvature bounds as in~\cite{bgrs14} and~\cite{IM14}, measures $\nu$ satisfying a (weaker) Poincar\'e inequality as in~\cite{FIL14}. Under these additional assumptions on $\nu$, the goal is then to obtain better constants in the logarithmic Sobolev inequality, mimicking in a sense the phenomenon observed in the Poincar\'e inequality when considering test functions orthogonal to the first eigenfunctions. In Indrei-Marcon \cite{IM14}, the deficit  is controlled  by the Wasserstein distance for the class of centered functions with upper and lower bounded curvature. The authors in~\cite{bgrs14} also give new bounds in terms of conditionally centered vectors. Further improvements are given in~\cite{FIL14} in terms of the $W_{1,1}$ distance defined in Section~\ref{tenso}. Here again our bounds share the advantages of holding without any smoothness, centering, etc. hypothesis on $\nu$, and of having the good dimensional behaviour when considering product measures.
 
\medskip
\section{Applications to Fokker-Planck equations}\label{sec-application-fp}

Let us now see how our results (or methods) lead to short-time smoothing of the entropy and improved contraction rates for the laws of solutions to stochastic differential equations.

\smallskip

For this, let again $V$ be a $C^2$ function on $\mathbb R^n$ such that $\int e^{-V}=1$ and $\He{V} \geq R\,{\rm Id}_n$, with $R$ possibly negative, and satisfying the doubling condition $V(x+y) \leq C(1 + V(x) + V(y))$ for a $C$ and all $x, y.$ Let also $\mu$ be the probability measure with density $e^{-V}$.
We let $u_0$ in $P_2(\R^n)$ and consider gradient flow solutions $u = (u_t)_{t \geq 0} \in C([0, + \infty), P_2(\R^n))$ of the Fokker-Planck equation
\begin{equation}\label{FP}
\frac{\partial u_t}{\partial t} = \Delta u_t + \nabla \cdot (u_t \nabla V), \qquad t>0, x \in \rr^n
\end{equation}
as in~\cite[Chap.~11.2.1]{ambrosio-gigli-savare} and \cite[Th.~4.20 and~4.21]{daneri-savare} (see also~\cite{lisini}). Equation~\eqref{FP} holds in the sense of distributions. Moreover,  by~\cite[Th.~11.2.8]{ambrosio-gigli-savare} or again~\cite{daneri-savare},  for any $t>0$ the solution $u_t$ has a density; for almost every $t>0$ this density is in $W^{1,1}_{loc} (\mathbb R^n)$, with $\nabla u_t/ u_t + \nabla V \in L^2(u_t)$; finally $t \mapsto I(u_t \vert \mu) \in L^1_{loc}(]0, + \infty[)$ and
$$
\frac{d}{dt} H (u_t \vert \mu) = - I(u_t \vert \mu)
$$
for almost every $t>0.$  The solution $u_t$ can be seen as the law at time $t$ of the solution $(X_t)_{t \geq 0}$ to the stochastic differential equation
$$
d X_t = \sqrt{2} \, dB_t - \nabla V (X_t) \, dt.
$$
Here $(B_t)_{t \geq 0}$ is a standard Brownian motion on $\mathbb R^n$ and the initial datum $X_0$ has law $u_0$.

Moreover, the interpretation of~\eqref{FP} as the gradient flow of $H(\cdot \vert \mu)$ on the space 
$P_2(\R^n)$ has enabled to obtain the following short-time and contraction properties (see~\cite[Th.~11.2.1]{ambrosio-gigli-savare} and \cite[Chap.~24]{villani-book1}). Let $u$ and $v$ be solutions to~\eqref{FP}. Then 
\begin{equation}\label{regH}
H(u_t \vert \mu) \leq \frac{W_2^2 (u_0, \mu)}{2t} e^{2 \max\{-R, 0\}\,t}, \qquad t > 0
\end{equation}
and
\begin{equation}\label{contr}
W_2(u_t, v_t) \leq e^{-Rt} \, W_2 (u_0, v_0), \qquad t \geq 0.
\end{equation}
In particular, if $R>0$, then $u_t$ converges to the steady state $\mu$ as
\begin{equation}\label{cvexpo}
W_2(u_t, \mu) \leq e^{-Rt} \, W_2 (u_0, \mu), \qquad t \geq 0.
\end{equation}

The purpose of this section is to improve these three properties by means of the tools and inequalities in the above sections.

 \subsection{Short-time smoothing of the entropy}\label{regu}
 
 In the Gaussian case where $\mu$ is the standard Gaussian measure $\gamma$, the solution to \eqref{FP} is given by the Mehler formula 
 (see~\cite[Sect.~2.7.1]{bgl-book}). In particular the fundamental solution, with initial datum $u_0$ the Dirac mass at $0$, is at time $t>0$ the Gaussian measure with variance $\sigma_t^2 = 1-e^{-2t}$:
 $$
 u_t(x) = (2 \pi \sigma_t^2)^{-n/2} e^{-x^2/(2 \sigma_t^2)}, \qquad z \in \mathbb R^n.
$$
Its relative entropy can be computed as 
$$
H(u_t \vert \gamma) = \int_{\mathbb R^n} u_t(x) \log \frac{u_t(x)}{\gamma(x)} dx = -\frac{n}{2} \big[ e^{-2t} + \log(1-e^{-2t}) \big].
$$
Of course this is coherent with \eqref{regH}, with $R=1$, since
$$
-\frac{n}{2} \big[ e^{-2t} + \log(1-e^{-2t}) \big] \leq \frac{n}{2t} = \frac{W_2^2 (u_0, \mu)}{2t}
$$
by direct computation. In fact, for $t \sim 0$ one can observe that 
$$
H(u_t \vert \gamma) \sim \frac{n}{2} \log \frac{1}{t} \cdot
$$

On the other hand, let $u$ be a solution to \eqref{FP}, still in the Gaussian case, and with initial datum $u_0$ such that $u_0(\vert x \vert^2) = n = \gamma (\vert x \vert^2).$ Then $u_t(\vert x \vert^2) = n$ for all $t$ since
 \begin{equation}\label{deuxmom}
\frac{d}{dt} \int \vert x \vert^2 \, du_t = 2n - 2 \, \int \vert x \vert^2 \, du_t.
\end{equation}
In particular, in the notation $H(t) = H(u_t \vert \gamma) / n$ and $I(t) = I(u_t \vert \gamma)/n,$ the dimensional Gaussian logarithmic
 Sobolev inequality~\eqref{eq-logsob-bakry-ledoux} simplifies as $ 2H \leq \log(1+I).$ Hence
$$
H'(t) = -I(t) \leq 1 - e^{2H(t)}, \qquad \textrm{for} \; a. e.\;  t >0.
$$
By the change of variable $x(t) = e^{-2h(t)}$ this integrates into
$$
x(t) e^{2t} \geq x(0) + e^{2t} - 1 \geq e^{2t} - 1.
$$
In other words
$$
 H(u_t \vert \gamma) \leq -\frac{n}{2} \log (1 - e^{-2t}), \qquad t>0
 $$
which gives the same short-time behaviour.
 
 \bigskip
 
 More generally :
 
 \begin{prop}\label{t=0}
Let $u$ be a solution to~\eqref{FP} with $\He{V} \geq R\,{\rm Id}_n, R>0$, and with initial condition $u_0$ in $P_2(\R^n).$ Let $T>0$ and assume that $u_t(\vert \nabla V \vert^2) \leq M$ for $t$ in $[0, T]$. Then there exists a constant $c>0$ depending only on $n, R$ and $M$ such that
 $$
 H(u_t \vert \mu) \leq \max \Big\{1, \frac{n}{2} \log \frac{c}{t} \Big\}, \qquad t \leq T.
 $$
 \end{prop}
\begin{remark}
The moment assumption $u_t(\vert \nabla V \vert^2) \leq M$ for $t$ in $[0, T]$, is not a restrictive condition. It can indeed be checked by time differentiating $u_t(|\nabla V|^2)$ and controlling its non explosion via a Lyapunov type condition on $u_0e^V$ or on derivatives of $V$ for instance. 

It can also be checked by observing that the Markov semigroup $(P_t)_{t\geq 0}$ with generator $L=\Delta-\nabla V\cdot \nabla$ is such that
$\int \phi du_t=\int P_t \phi du_0$
 for any test function $\phi$. In particular, if $\Phi$ is a convex function and if the initial datum has a density also denoted $u_0$, then
\begin{multline*}
u_t(\vert \nabla V \vert^2)=\int \vert \nabla V \vert^2 du_t=\int P_t(\vert \nabla V \vert^2)du_0=\int P_t(\vert \nabla V \vert^2)u_0e^Vd\mu\\
\leq\int \Phi(P_t(\vert \nabla V \vert^2))d\mu+\int \Phi^*(u_0e^V)d\mu\leq \int \Phi(\vert \nabla V \vert^2)d\mu+\int \Phi^*(u_0e^V)d\mu.
\end{multline*}
Here we use the fact that $t\mapsto \int \Phi(P_t(\vert \nabla V \vert^2))d\mu$ is non increasing since $\Phi$ is convex. The moment assumption is then satisfied for all $T>0$ as soon as the right hand side is finite for a convex function $\Phi$. 
\end{remark}
 
\begin{eproof}
We shall let $c$ denote diverse positive constants depending only on $n$, $M$ and $R$. By Corollary~\ref{coro-lsbe} applied to the measure $f^2 \mu = u_t$,  and integration by parts, there holds
$$
H(u_t \vert \mu) \leq n(s-1-\log s) + \frac{1-s^2}{2R} u_t (\vert \nabla V \vert^2) + \frac{s(s-1)}{R} u_t ( \Delta V) + \frac{s^2}{2R} I(u_t \vert \mu)
$$
for $t > 0$ and $s>0$. Recall that $I$ has been introduced in~\eqref{eq-def-fisher}. Since $V$ is convex, then $\Delta V \geq 0$ and then
$$
H(u_t \vert \mu) \leq - n \, \log s + c + \frac{s^2}{2R} I(u_t \vert \mu)
$$
for all $s \in ]0, 1]$ and $t \in ]0, T].$

Now, as far as $H(t):= H(u_t \vert \mu) \geq 1$, then $I (t):= I(u_t \vert \mu) \geq 2R$ so that $s = \sqrt{2R/I}$ is smaller than $1$. For this $s$ we obtain
$$
H \leq c + \frac{n}{2} \log I.
$$
Hence
$$
H'(t) = - I (t) \leq -e^{2H(t)/n - c}
$$
 for almost every $t>0.$ As above $x(t) = e^{-2H/n}$ satisfies $x(t) \geq x(0) + ct \geq ct$ by time integration. Written in terms of~$H$, this concludes the proof.
\end{eproof}

\subsection{Refined contraction properties}


Let us now see how to make~\eqref{contr} finer. Still by~\cite[Th.~8. 3. 1]{ambrosio-gigli-savare} and \cite[Th.~4.20 and~4.21]{daneri-savare}, one can write \eqref{FP} as the continuity equation
$$
\frac{\partial u_t}{\partial t} + \nabla \cdot (\xi[u_t] u_t) =0,  \qquad t>0, x \in \rr^n
$$
with $\xi[u_t] = - \nabla V - \nabla \log u_t$. Then for almost every $t>0$
\begin{eqnarray}
\! \! \! \! \! \! \! \! \! \! \! - \frac{1}{2} \frac{d}{dt} W_2^2(u_t, v_t) 
\! \! \! \! 
&=&  
 \int \big( \xi[v_t](\nabla \varphi_t(x)) - \xi[u_t](x) \big) \cdot (\nabla \varphi_t(x) - x ) \, u_t(x) \, dx\label{egal}
 \\
 &\geq&
 \! \! \! \! \! \!
 \int
\! \! \Big[ \Delta \varphi_t (x) \!+\! \Delta  \varphi_t^*(\nabla \varphi_t(x)) \! - \! 2n + \! \big( \nabla V (\nabla \varphi_t(x)) \!-\! \nabla V (x) \big) \! \cdot \! (\nabla \varphi_t(x) \!-\!  x) \Big] \, u_t(x) \, dx\label{inegal}
\end{eqnarray}
for two solutions $u$ and $v.$ Here $\varphi_t$ is the  convex map such that $v_t = \nabla \varphi_t \# u_t$ and $u_t = \nabla \varphi_t^* \# v_t$ for the Legendre transform  $\varphi_t^*$ of $ \varphi_t$ (see Section~\ref{rappels}). Equality~\eqref{egal} follows from~\cite[Th.~23.9]{villani-book1} (see also~\cite[Th.~8.4.7]{ambrosio-gigli-savare}); its assumptions are satisfied since (and likewise for $v$)
$$
\int_{t_1}^{t_2} \int_{\mathbb R^n} \vert \xi[u_s] \vert^2 du_s \, ds = \int_{t_1}^{t_2} I(u_s \vert \mu) ds =  H(u_{t_1} \vert \mu) - H(u_{t_2} \vert \mu) \leq H(u_{t_1} \vert \mu)
$$
which is finite for any $t_2 > t_1 >0, $ as observed above. Inequality~\eqref{inegal} follows from a weak integration by parts, as in \cite[Th. 1.5]{lisini}; there again $\Delta \varphi_t$ is the trace of the Alexandrov Hessian of $\varphi_t$.

 Now, for given $t>0$ and $u_t$-almost every $x$, the symmetric matrix $\He {\varphi_t} (x)$ is positive, as recalled in Section~\ref{rappels} : letting $e^{2 \lambda_i(x)}$ for $i= 1, \dots, n$ its $n$ positive eigenvalues , then its inverse matrix $\He{\varphi_t^*}(\nabla \varphi_t(x))$ (see again~Section~\ref{rappels}) has eigenvalues $e^{-2 \lambda_i(x)}$; hence at point $x$
\begin{equation}\label{eq-laplace}
\Delta \varphi_t+ \Delta  \varphi_t^*(\nabla \varphi_t) - 2n
= \tr \big[ \He{ \varphi_t}  \big]+ \tr \big[ \He{ \varphi_t^*}(\nabla \varphi_t) \big] - 2n = \sum_i \big( e^{2 \lambda_i} + e^{-2 \lambda_i} - 2) = 4 \sum_i \sinh^2(\lambda_i).
\end{equation}
Hence, by convexity of $\sinh^2$ and the Jensen inequality, and \eqref{entdx},
\begin{eqnarray*}
\int \Big[ \Delta \varphi_t (x)+ \Delta  \varphi_t^*(\nabla \varphi_t(x)) - 2 \, n \Big] u_t(x) \, dx
&=& 
4 n \, \frac{1}{n} \sum_i \int \sinh^2(\lambda_i(x)) \, u_t(x) \, dx 
\\
&\geq &
4n \sinh^2 \left( \frac{1}{n} \sum_i \int \lambda_i (x) \, u_t(x) \, dx \right)
\\
&=&
 4 n \sinh^2 \left( \frac{1}{2n} \int \log \det \He{\varphi_t}(x) \, u_t(x) \, dx \right)
\\
&=& 4n \sinh^2 \left( \frac{\entf{dx} (v_t) - \entf{dx}(u_t)}{2n} \right).
\end{eqnarray*}
Since $\He{V} \geq R\,{\rm Id}_n,$ we obtain
\begin{equation}\label{derW22}
- \frac{1}{2} \frac{d}{dt} W_2^2(u_t, v_t) \geq 4n \sinh^2 \left( \frac{\entf{dx} (v_t) - \entf{dx}(u_t)}{2n} \right)+ R W_2^2(u_t, v_t).
\end{equation}
By time integration this ensures the following dimensional contraction property :
\begin{prop}\label{propcontr}
In the above notation, if $\He{V} \geq R\,{\rm Id}_n$ for $R\in\R,$ then for any solutions to \eqref{FP}
\begin{equation}\label{cvw2}
W_2^2(u_t, v_t) \leq e^{-2Rt} W_2^2(u_0, v_0)  - 8n \int_0^t e^{-2R(t-s)}\sinh^2 \left( \frac{\entf{dx} (v_s) - \entf{dx}(u_s)}{2n} \right)  \, ds, \qquad t \geq 0.
\end{equation}
\end{prop}

 For the heat equation, namely for $V=0$, then the associated Markov generator $L = \Delta$ satisfies the $CD(0,n)$ curvature-dimension condition: in particular in this case the bound~\eqref{cvw2} has been derived in~\cite{BGG13} and~\cite{BGGK}, and is also a consequence of~\cite{EKS13}. For $V \neq 0$, then the associated generator $L = \Delta - \nabla V \cdot \nabla$ satisfies a $CD(R, \infty)$ but no $CD(R,n)$ condition: in particular the bound~\eqref{cvw2} can not be obtained from the works mentioned above.

\begin{remark}
The above computation can be extended to drifts $A(x)$ which are not gradients. In this case the assumption $\He{V} \geq R\,{\rm Id}_n$ should be replaced by the monotonicity condition
$(A(y) - A(x) ) \cdot (y-x) \geq R \, \vert y-x \vert^2$ for all $x,y$ (see~\cite{BGG11} for this non-gradient case). 
\end{remark}

\subsection{A formal gradient flow argument to Proposition \ref{propcontr}}

In this subsection, we provide an alternative formal argument to Proposition \ref{propcontr} based on gradient flow.
\medskip

We begin with the following elementary lemma which gives additional information to \cite[Lem.~2.2]{EKS13}. 

\begin{lemma}\label{lemmefonda} Let $\psi$ be a $C^2$ function on $[0,1]$. Then the following properties are equivalent:
\begin{itemize}
\item  $\psi'' \geq \psi'^2/n$;
\item  for all $r, s$ in $[0,1]$,
\begin{equation}\label{lemme1}
n - \psi'( r) (s-r) \geq n \, e^{\frac{\psi( r)- \psi( s)}{n}};
\end{equation}
\item  for all $r, s$ in $[0,1]$,
\begin{equation}\label{lemme2}
\big( \psi'(s) - \psi'( r) \big) (s-r) \geq 4n \sinh^2 \Big( \frac{\psi(s)- \psi( r)}{2n} \Big).
\end{equation}
\end{itemize}
\end{lemma}

\begin{eproof}
Let indeed $U= e^{-\psi/n}$, so that 
$$
U'' = - \Big( \psi'' - \frac{\psi'^2}{n} \Big) \frac{U}{n}.
$$ 
Then $\psi'' \geq \psi'^2/n$ if and only if $U$ is concave, hence if and only 
$$
e^{-\frac{\psi(s)}{n}} = U(s) \leq U( r) + U'( r) (s-r) = e^{-\frac{\psi( r)}{n}} - \frac{\psi'( r)}{n} e^{-\frac{\psi( r)}{n}} (s-r)
$$
for all $r, s \in [0,1]$, which is  \eqref{lemme1} when multiplying both sides by $e^{\psi(r)/n}$.

Adding \eqref{lemme1} with the corresponding bound obtained with $r, s$ instead of $s, r$ leads to \eqref{lemme2}.
Conversely, dividing \eqref{lemme2} by $(s-r)^2$ and letting $s$ go to $r$ gives $\psi'' \geq \psi'^2/n$ at point $r$.
\end{eproof}

\bigskip
 
Let now $\mu^0$ and $\mu^1$ be absolutely continuous measures in $P_2(\mathbb R^n)$, $\nabla \varphi$ their Brenier map  and $(\mu^s)_{s \in [0, 1]}$ the geodesic between them, as in Section~\ref{rappels}. Here again  we identify the measures with their densities.
Let us now recall why the function $\psi : s \mapsto \entf{dx} (\mu^s)$ formally satisfies $\psi'' \geq \psi'^2/n$  on [0,1].
For this, recall from Section~\ref{rappels} that for $\mu^0$-almost every $x$ the Alexandrov Hessian $\He{\varphi} (x)$ is positive, so that the eigenvalues $\theta_i(x)$ of $\He{\varphi} (x) -I$ are $> -1.$
Writing \eqref{entdx} with the measures $\mu_0 = \mu^0$ and $\mu_1 = \mu^s,$ we obtain
$$
\psi(0) 
= \psi(s) + \int \log \det(I + s(\He{\varphi} (x)-I)) \, d\mu^0 (x)
= \psi (s) + \sum_i \int \log (1+s \theta_i(x)) \, d\mu^0 (x).
$$
Hence
\begin{equation}\label{psi'}
\psi'(s) = - \sum_i \int \frac{\theta_i}{1 + s \theta_i} d\mu^0
\end{equation}
and then by the Cauchy-Schwarz inequality
$$
\psi''(s) = n \, \frac{1}{n} \sum_i \int \frac{\theta_i^2}{(1 + s \theta_i)^2} d\mu^0 
\geq n \left( \frac{1}{n}  \sum_i \int \frac{\theta_i}{1 + s \theta_i} d\mu^0 \right)^2 
 = \frac{1}{n} \psi'(s)^2.
$$

\begin{remark}\label{rem33}
Identity \eqref{psi'} can also be formally checked using the continuity equation solved by $(\mu^s)_{s \in [0,1]}$:
$$
\frac{\partial \mu^s}{\partial s} + \nabla \cdot (\mu^s v^s) =0.
$$
 Here the vector field $v^s$ satisfies $v^s(x + s(\nabla \varphi(x) -x)) = \nabla \varphi(x) - x$, see e.~g.~\cite[Th.~5.51]{villani-otp}. For, and recalling that $\psi(s)= \int \mu^s \log \mu^s \, dx$
$$
\psi'(s) 
= \! - \! \int \nabla \cdot (v ^s \mu^s) \log \mu^s \, dx 
= \! - \! \int \nabla \cdot v^s \, \mu^s \, dx\\
= \! - \! \int \big( \nabla \cdot v^s \big) (x + s (\nabla \varphi(x)-x)) \, d\mu^0 (x)
$$
by integration by parts and since $(x + s (\nabla \varphi(x)-x)) \# \mu^0 = \mu^s$. Identity \eqref{psi'} follows since by chain rule
$$
\big( \nabla \cdot v^s \big) (x + s (\nabla \varphi(x)-x)) = \tr \Big[ (\He{\varphi} (x) - I) \big(I + s (\He{\varphi} (x)-I) \big)^{-1} \Big] = \sum_i \frac{\theta_i}{1 + s \theta_i} \cdot
$$
\end{remark}

\begin{remark} In the above notation, observe that \eqref{lemme1} in Lemma \ref{lemmefonda} for $\psi(s) = Ent_{dx}(\mu^s), r=0$ and $s=1$ formally leads to \eqref{bornelaplace} in Lemma \ref{lemma-laplace}. For, in the notation of Remark~\ref{rem33} and by integration by~parts,
$$
\psi'(0) = \int \nabla \mu^0 \cdot v^0 \, dx  = \int \nabla \mu^0 \cdot (\nabla \varphi -x) dx= n - \int \Delta \varphi \, d\mu^0.
$$
\end{remark}

\bigskip

We can now deduce an alternative formal argument to the bound in Proposition \ref{propcontr}.

\smallskip

We begin with the following classical observation in Euclidean space : Let $X$ and $Y$ be two solutions of the Euclidean gradient flow $X'_t = - \nabla U(X_t)$ in $\mathbb R^d$, where $U : \mathbb R^d \to \mathbb R$ is a smooth potential. For $t >0$ let $U_t (s) = U(X_t + s(Y_t - X_t))$ for $s \in [0, 1]$. Then
\begin{equation}\label{ddt}
 - \frac{1}{2} \frac{d}{dt} \vert Y_t - X_t \vert^2
=
(Y_t - X_t) \cdot (\nabla U(Y_t) - \nabla U(X_t))
=
U_t' (1) - U_t'(0).
\end{equation}

Let now $u$ and $v$ two solutions to the Fokker-Planck equation~\eqref{FP}, which by~\cite[Chap.~11.2]{ambrosio-gigli-savare} and~\cite[Chap.~23]{villani-book1} is the gradient flow of $H(\cdot \vert \mu)$ on the space 
$P_2(\R^n)$. For any $t > 0$, let $\nabla \varphi_t$ be the optimal transport map between $u_t$ and $v_t$, and $(\mu_t^s)_{s \in [0,1]}$ be the geodesic path in $P_2(\R^n)$ 
between $u_t$ and $v_t$, as in Section~\ref{rappels}. 
Then, formally and by analogy with~\eqref{ddt},
\begin{equation}\label{ddtE}
- \frac{1}{2} \frac{d}{dt} W_2^2(u_t, v_t)
=
E_t' (1) - E_t'(0)
\end{equation}
where for given $t>0$ we let 
$$
E_t (s) = H(\mu_t^s \vert \mu) = \entf{dx} (\mu_t^s) + \int V \, d\mu_{t}^s.
$$
Indeed, let $\psi : s \mapsto \entf{dx} (\mu_t^s)$ for given $t$ and, for each $x$ let the matrix $\He{\varphi_t} (x)$ have eigenvalues $e^{2\lambda_i}$.  Then, in the above notation $\theta_i = e^{2\lambda_i} -1$,~\eqref{psi'} for $\mu = u_t$ gives
\begin{eqnarray*}
\psi'(1) - \psi'(0) 
&=&
\! \int \sum_i \Big[ \theta_i - \frac{\theta_i}{1+\theta_i} \Big] du_t 
\!\\
&=&
\! \int \sum_i \frac{\theta_i^2}{1+\theta_i} du_t 
 \\
&=&
\! \int \sum_i \Big[ e^{2 \lambda_i} + e^{-2\lambda_i} - 2  \Big] du_t
= 
\! \int \Big[ \Delta \varphi_t + \Delta  \varphi_t^*(\nabla \varphi_t) - 2 \, n \Big] du_t
\end{eqnarray*}
as in~\eqref{eq-laplace}. Using moreover the formal derivative
$$
\frac{d}{ds} \int V \, d\mu_{t}^s = \frac{d}{ds} \int V \big( x + s(\nabla \varphi_t(x) - x) \big) du_t(x) = \int  \nabla V \big( x + s(\nabla \varphi_t(x) - x) \big) \cdot (\nabla \varphi_t(x) - x) \, du_t(x)
$$
for $s=0, 1$ we formally recover~\eqref{inegal} in~\eqref{ddtE}. 

We now use the fact that for given $t$ the function $\psi$ satisfies $\psi'' \geq \psi'^2/n$  on $[0,1]$. Then, by \eqref{lemme2} in Lemma \ref{lemmefonda} for $r=0$ and $s=1$ we obtain
$$
E_t'(1) - E_t'(0) \geq 4n \sinh^2 \Big(\frac{Ent_{dx}(v_t) - Ent_{dx} (u_t)}{2n} \Big) + \int \big( \nabla V (\nabla \varphi_t(x)) - \nabla V(x) \big) \cdot \big( \nabla \varphi_t(x) - x \big) \, du_t (x).
$$
Since $\int \vert \nabla \varphi_t(x) - x \vert^2 du_t(x) = W_2^2(u_t, v_t)$ this leads to~\eqref{derW22} and then to~\eqref{cvw2} as soon as $\He{V} \geq R\,{\rm Id}_n.$

\bigskip


\subsection{Improved convergence rates}

In this section we consider a solution $u$ to~\eqref{FP} in the Gaussian case where $\mu =\gamma,$ and for which we can take $R=1$ above. Let us see how the contraction property~\eqref{cvw2} can make the convergence estimate~\eqref{contr} more precise. 

\smallskip

For simplicity we assume that $u_0(\vert x \vert^2) \leq n = \gamma( \vert x \vert^2).$ Then  $u_t(\vert x \vert^2)  \leq n$ for all $t$, by~\eqref{deuxmom}.
Hence \eqref{entropies} and the Talagrand inequality \eqref{tala} ensure that $0 \leq W_2^2(u_t, \gamma) < 2n$ and 
$$
\frac{\entf{dx} (u_t) - \entf{dx} (\gamma)}{n} \geq - \log \Big( 1 - \frac{W_2^2(u_t, \gamma)}{2n} \Big).
$$
In particular the right-hand side is non negative. 
Moreover, for the stationary solution $v_t = v_0 = \gamma$, the contraction property \eqref{cvw2} with $R=1$, in the form \eqref{derW22}, implies
$$
- x' \geq   \frac{x^2}{1-x} + 2 \, x
$$
where $x(t) = W_2^2(u_t, \gamma) / (2n) \in [0,1)$. Here we use that $\sinh (\log x) = (x - 1/x)/2.$
In other words $z(t) = 1 - (1-x(t))^2$ satisfies $z' \leq - 2z$. This integrates into $z(t) \leq e^{-2t} z(0)$, that is,
%
\begin{equation}\label{eqx}
x(t) \leq 1 - \Big(1 - (2 x(0) - x(0)^2) e^{-2t} \Big)^{\frac12}.
\end{equation}
By the lower bound
\begin{equation}\label{racine}
1 - (2 x(0) - x(0)^2) e^{-2t} \geq (1-x(0)e^{-2t})^2
\end{equation}
it implies the classical bound~\eqref{contr}. It also improves it:
for instance~\eqref{eqx} can be written as
$$
W_2^2(u_t, \gamma) 
\leq  W_2^2(u_0, \gamma)  e^{-2t} \frac{2-x(0)}{1 + \Big(1 - (2 x(0) - x(0)^2) e^{-2t} \Big)^{\frac12}}. 
$$
Then by~\eqref{racine} we obtain 
\begin{coro}\label{coro-cvOu}
In the above notation, let $u$ be a solution to~\eqref{FP} in the Gaussian case, with initial datum $u_0$ such that $u_0 (\vert x \vert^2) \leq n.$ Then for all $t \geq 0$
$$
W_2^2(u_t, \gamma) 
\leq W_2^2(u_0, \gamma) e^{-2t} \, \frac{1-W_2^2(u_0, \gamma)/(4n)}{1 - W_2^2(u_0, \gamma) e^{-2t}/(4n)}.
$$
\end{coro}
Observe that the quotient is smaller than $1$.

\medskip

\begin{remark}
The Gaussian assumption is used here only to ensure uniform convexity of the potential (hence the Talagrand inequality), and that $\int V du_t \leq \int V e^{-V}dx$ as soon as this holds at $t=0$.
\end{remark}

\section{Brascamp-Lieb inequalities}
\label{sec-dim-BL-ine}

It is classical that linearizing a logarithmic Sobolev inequality leads to a Poincar\'e inequality, which in the Gaussian case is the Brascamp-Lieb inequality. In this section we shall see how to obtain two different dimensional Brascamp-Lieb inequalities:  a first one by an improvement of the classical $L^2$ method, and a second one by linearization in the Borell-Brascamp-Lieb inequality~\eqref{eq-borell-bl}. 

\subsection{Brascamp-Lieb inequality by $L^2$ method}

\begin{prop}[Dimensional Brascamp-Lieb inequality I]
\label{prop-BL}
Let $\mu$ be a probability measure on $\rr^n$ with density $e^{-V}$ where $V$ is a $\mathcal C^2$ function satisfying $\He{V}>0$.
Then
\begin{equation}\label{BLn}
\varf{\mu}(f) \leq \int \nabla f \cdot \He{V}^{-1} \nabla f \, d\mu - \frac {\Big( \int V \, f \, d\mu - \int V \, d\mu \int f \, d\mu\Big)^2}{n - \varf{\mu}(V)}
\end{equation}
for all $\mathcal C^1$ compactly supported functions $f$.
\end{prop}


\begin{remark}
V. H. Nguyen~\cite{nguyen} has proven that  $\varf{\mu} (V) \leq n$ for $V$ convex. We will observe in the proof that even $Var_{\mu} (V) < n$ as soon as $\He{V}  >0$. In fact, it follows from the bound~\eqref{BLn} for $f=V$ that $Var_{\mu} (V) \leq \frac{nI}{n+I} < n$ where $I =  \int \nabla V \cdot \He{V}^{-1} \nabla V \, d\mu.$ In particular, if $R\,{\rm Id}_n \leq \He{V} \leq S\,{\rm Id}_n,$ then $I \leq R^{-1} \int \vert \nabla V \vert^2 d\mu =
R^{-1} \int \Delta V d\mu \leq nS/R$ and $\varf{\mu} (V) \leq \frac{nS}{R+S}.$ The latter inequality is an equality (to $n/2$) for the Gaussian measure with any variance, for which $R=S$.
\end{remark}

If $\mu=\gamma$ is the standard Gaussian measure then \eqref{BLn} is exactly the dimensional (Poincar\'e) inequality~\eqref{Pdimgauss} (and in particular equality holds for $f=\vert x \vert^2/2$). 

\medskip

In the non Gaussian case, G. Harg\'e has derived the following improvement of the Brascamp-Lieb inequality, see~\cite[Th.~1]{harge-jfa} : if $V$ is a  $\mathcal C^2$ function satisfying  $R\,{\rm Id}_n \leq \He{V}  \leq S\,{\rm Id}_n$ for constants $0 \leq R \leq S$, then
\begin{equation}\label{harge}
\varf{\mu}(f) \leq \int \nabla f \cdot \He{V}^{-1} \nabla f \, d\mu - \frac{1+R/S}{n} \Big( \int V \, f \, d\mu - \int V \, d\mu \int f \, d\mu \Big)^2
\end{equation}
for all  $f$.

 We do not know in full generality which of the coefficients $(n- Var_{\mu}(V))^{-1}$ and  $n^{-1} (1+R/S)$ in the corrective terms of  \eqref{BLn} and \eqref{harge} is the larger.
 
Besides being equal (to $2n^{-1}$) in the Gaussian case, both coefficients are always larger than $n^{-1}$. More precisely the coefficient in \eqref{BLn} is always strictly larger than $n^{-1}$ whereas the coefficient in \eqref{harge} is $n^{-1}$ when $R=0$ (no uniform convexity) or $S = + \infty$ (no upper bound on $\He{V}$): hence at least in these cases our bound is stronger. 

The bound~\eqref{harge} has been obtained in~\cite{harge-jfa} by a $L^2$ argument. We shall see in the appendix that it can be formally recovered by linearization in the Monge-Amp\`ere equation.

\bigskip

\noindent
{\bf Proof of Proposition \ref{prop-BL}.}  Let $\omega$ be in the space $\mathcal C^{\infty}_c$ of $\mathcal C^{\infty}$ and compactly supported functions. Then
$$
\int ||\He{\omega}||^2_{HS} \, d\mu-\frac{1}{n}\PAR{\int \Delta \omega \, d\mu}^2\geq 0
$$  
by the Cauchy-Schwarz inequality; here $||\He{\omega}||^2_{HS}=\sum_{i,j=1}^n(\partial_{ij}\omega)^2$ is the squared Hilbert-Schmidt norm of the matrix $\He{\omega} = (\partial_{ij}\omega)_{i, j}$. In other words
\begin{equation}
\label{eq-l2}
\int \nabla \omega\cdot\He{V}\nabla \omega \, d\mu\leq \int \PAR{||\He{\omega}||^2_{HS}+\nabla \omega\cdot\He{V}\nabla \omega} \, d\mu-\frac{1}{n}\PAR{\int \Delta \omega \, d\mu}^2.
\end{equation}
Moreover, by integration by parts, 
\begin{equation}\label{gamma2}
\int \PAR{||\He{\omega}||^2_{HS}+\nabla \omega\cdot\He{V}\nabla \omega} \, d\mu=\int (L\omega)^2d\mu, \qquad \int \Delta \omega \, d\mu=-\int VL\omega \, d\mu
\end{equation}
with $L=\Delta-\nabla V\cdot\nabla$, see~\cite[Section 3.2]{bgl-book}.

Let now $f$ be a $\mathcal C^1$ compactly supported function. Then pointwise 
$$
 2 \, \nabla f\cdot \nabla \omega\leq \nabla \omega\cdot\He{V}\nabla \omega+\nabla f\cdot \He{V}^{-1}\nabla f. 
$$ 
From these remarks, inequality~\eqref{eq-l2} implies  
$$
2\int \nabla f \cdot \nabla \omega \, d\mu \leq \int \nabla  f\cdot\He{V}^{-1} \nabla f \, d\mu +  \int (L\omega)^2d\mu-\frac{1}{n}\PAR{\int  VL\omega \, d\mu}^2.
$$ 
Let now $h = - L\omega$. Then $\int \nabla f\nabla \omega \, d\mu=-\int fL \omega \, d\mu=\int fh \, d\mu$ by integration by parts. 

To sum up, we have obtained
\begin{equation}
\label{eq-l21}
2\int fh \, d\mu \leq \int \nabla  f\cdot\He{V}^{-1} \nabla f \, d\mu +  \int h^2d\mu-\frac{1}{n}\PAR{\int  Vh \, d\mu}^2
\end{equation}
for any $h$ in $L(\mathcal C^{\infty}_c)$ and any $\mathcal C^1$ compactly supported function $f$. 

But, by~\cite[Lem.~9]{harge-jfa} for instance, $L(\mathcal C^{\infty}_c)$ is dense (for the $L^2(\mu)$ norm) in the space of functions $h\in L^2(\mu)$ such that  $\int hd\mu=0$. Hence, formula~\eqref{eq-l21} extends to any $h\in L^2(\mu)$ such that  $\int h \, d\mu=0$.

In particular, given a $\mathcal C^1$ compactly supported function $f$ such that $\int f \, d\mu=0,$ we can apply~\eqref{eq-l21} to $ h = f + a \big(V - \int V \, d\mu \big)$ with $a\in\R$. Observe indeed that $V \in L^2(\mu)$ for $\mu = e^{-V}$ with $V$ convex. We get
$$
\int f^2 \, d\mu \leq  \int  \nabla f \cdot \He{V}^{-1} \nabla f \, d\mu + I \, a^2 - 2 a \; \frac{\varf{\mu} (V)}{n}  \int V f \, d\mu - \frac{1}{n} \Big( \int V f \, d\mu \Big)^2
$$
for all $a$, where $\displaystyle I = \varf{\mu} (V) (n - \varf{\mu} (V)) /n $. Necessarily $I$ is positive, that is, $\varf{\mu} (V) < n$. Indeed, if $I$ was non positive, then the left-hand side would be $-\infty$ by letting $a$ tend to $\pm \infty$, which is impossible. 

We finally optimise over $a$, choosing $\displaystyle a =  \int V f \, d\mu / (n - Var_{\mu} (V))$. This concludes the proof of Proposition~\ref{prop-BL} for any $f$ such that $\int f \, d\mu=0,$ and then for any $f$. 
\proofend


%

\subsection{Brascamp-Lieb inequality via the Borell-Brascamp-Lieb inequality}
\label{sec-dim-bb-ine}

The following result gives an improved version of the Brascamp-Lieb inequality~\eqref{BL} from the Borell-Brascamp-Lieb inequality. 
\begin{thm}[Dimensional Brascamp-Lieb inequality II]
\label{thm-brascamp-lieb-2}
Let $\mu$ be a probability measure on $\rr^n$ with density $e^{-V}$ where $V$ is a $\mathcal C^2$ function satisfying $\He{V}>0$. 
Then for any $\mathcal C^1$ and compactly supported function $f$ such that $\int fd\mu=0$, 
\begin{equation}\label{BLn-2}
\varf{\mu}(f) \leq \int \nabla f \cdot \He{V}^{-1} \, \nabla f \, d\mu - 
\int \frac{(f-\nabla f \cdot\He{V}^{-1} \, \nabla V)^2}{n+{\nabla V \cdot \He{V}^{-1} \, \nabla V}}d\mu.
\end{equation}
\end{thm}

\medskip

Theorem \ref{thm-brascamp-lieb-2} is proved in Appendix \ref{app-BL2}.

\medskip

For the standard Gaussian measure,  we obtain
\begin{coro}
The Gaussian measure $\gamma$ satisfies the dimensional Poincar\'e inequality
\begin{equation}\label{BLn-2-gaussian}
\varf{\gamma}(f) \leq \int |\nabla f |^2 \, d\gamma- 
\int \frac{(f-\nabla f\cdot x)^2}{n+|x|^2}d\gamma
\end{equation}
for any $\mathcal C^1$ and compactly supported function $f$ such that $\int fd\gamma=0$. 
\end{coro}
By the Cauchy-Schwarz inequality and integration by part, 
$$
\int \frac{(f-\nabla f\cdot x)^2}{n+|x|^2}d\gamma\geq \frac{\PAR{\int \nabla f\cdot xd\gamma}^2}{2n}=\frac{\PAR{\int \Delta fd\gamma}^2}{2n}=\frac{\PAR{\int f|x|^2/2d\gamma}^2}{n-Var_\gamma(|x|^2/2)} \cdot
$$
Therefore, for the Gaussian measure, inequality~\eqref{BLn-2-gaussian} is stronger than~\eqref{Pdimgauss} mentionned in the introduction (and naturally equality still holds for $f=|x|^2/2$).

\medskip

\subsection{Comparison of Brascamp-Lieb inequalities}
Many dimensional Brascamp-Lieb inequalities have recently been proved, and should be compared. We have already compared our inequality~\eqref{BLn} with G. Harg\'e's bound, as the same covariance term appears. Let us now compare \eqref{BLn-2} with other inequalities. It seems difficult to obtain a global comparison and we are only able to give partial answers or hints. 

\begin{itemize}

\item The present paper proposes the two inequalities~\eqref{BLn} and~\eqref{BLn-2}. In the Gaussian case we have already observed that~\eqref{BLn-2}-\eqref{BLn-2-gaussian} is stronger than~\eqref{BLn}. A variant of this argument shows that it is also the case for instance when $V(x)=x^{2a}+\beta, x \in \mathbb R$ with $a\in\N^*$ and  a normalisation constant $\beta$. We believe that it is the case for any $V$ since the additional term in~\eqref{BLn} vanishes for functions $f$ for which the one in~\eqref{BLn-2} does not. 

In fact, for a $\mathcal C^1$ function $f$ such $\int f e^{-V} = 0$, the additional term in~\eqref{BLn-2} vanishes if and only if there exists $a \in \mathbb R^n$ such that $f = a \cdot \nabla V  $ (and then $a = \int f (x) x e^{-V(x)}).$ For, if $f = \nabla f \, \He{V}^{-1} \nabla V$ on $\mathbb R^n$, then $g(y) = f(\nabla V^*(y))$ solves $g (y)= \nabla g (y) \cdot y$ on $\mathbb R^n$. Hence for fixed $y \in \mathbb R^n$ the map $t \mapsto g(ty) /t$ is constant; for $t=1$ and $t\to0$ this implies $g(y) = \nabla g(0) \cdot y$. This finally gives $f$, and conversely. But it is classical that these functions $f$ are exactly those for which equality holds in the Brascamp-Lieb inequality~\eqref{BL}. Hence the additional term in~\eqref{BLn-2} can be seen as a (weighted) way of measuring the distance of a function to the optimisers in the Brascamp-Lieb inequality~\eqref{BL}. 

Very recently, and under the same hypothesis as in Theorem~\ref{thm-brascamp-lieb-2},  D. Cordero-Erausquin  in~\cite[Prop.~6]{cordero15} proved that
\begin{equation}\label{dario}
Var_{\mu}(f) \leq \int \nabla f \cdot \He{V}^{-1} \, \nabla f \, d\mu - 
c \lambda(\mu) \int \He{V}^{-1} (\He{V} + c \lambda(\mu) {\rm Id}_n) ^{-1} \nabla f_0  \cdot \nabla f_0 d\mu
\end{equation}
for all $f$ satisfying $\int fd\mu=0$; here $f_0 = f - \int y f(y) d\mu(y) \cdot \nabla V$, $c$ is a numerical constant and 
 $\lambda(\mu)$ is the Poincar\'e constant of the measure $\mu$.  The additional term in~\eqref{dario} vanishes if and only if $f_0$ is a constant, so also appears here as a distance to the optimisers. A quantitative comparison between~\eqref{BLn-2} and~\eqref{dario} can not easily be performed as in particular a numerical constant appears in~\eqref{dario}. After the present work was completed, M.~Arnaudon, M. Bonnefont and A. Joulin~\cite{ABJ16} have derived Brascamp-Lieb inequalities in which the energy has been modified, instead of keeping the original energy and allowing for a remainder term, as here. We could not compare their results with ours.

\item We now turn to the Gaussian case when $\mu = \gamma$. We have already observed that~\eqref{BLn-2} is stronger  that~\eqref{BLn}, which is exactly~\eqref{Pdimgauss}. On the other hand,~\eqref{NPdimgauss} is a purely spectral inequality. We have numerically checked that~\eqref{BLn-2} 
implies~\eqref{NPdimgauss} for the Hermite polynomial functions $H_k$, $k\in\{1,\cdots, 7\}$. We believe that it is the case for all functions, but we do not have a proof of it. 

Let us conclude by mentioning the inequality
$$
Var_{\gamma}(f) \leq 6\int |\nabla f |^2 d\gamma - 
6\int \frac{(\nabla f\cdot x)^2}{n+|x|^2}d\gamma. 
$$
has been proved in~\cite[Sect.~2]{BL09}. Their extremal functions have been lost  since there is no equality when $f(x)=a\cdot x$ and the constant in front of the energy is larger than in our bounds. 


\end{itemize}

\begin{appendix}
 
\section{Proof of Theorem~\ref{thm-super-etoile}}
\label{sec-appendix-1}

{\bf Optimality of inequality~\eqref{eq-wetoile}.} When $W$ is strictly convex and satisfies 
(H1),
then~\eqref{legendre} holds, so that 
$$
\int \frac{W^*(\nabla W)}{W^{n+1}} \, dx=\int \frac{\nabla W\cdot x-W}{W^{n+1}} \, dx=-\frac{1}{n}\int \nabla (W^{-n})\cdot x \, dx-1 =0.
$$
In the last equality, we used an integration by parts, valid from hypothesis (H1) satisfied by $W$. This gives the equality case in~\eqref{eq-wetoile} when $g=W.$

\medskip

{\bf Proof of inequality~\eqref{eq-wetoile}.} Globally, the proof follows~\cite{bobkov-ledoux-sob}, but for completeness we give its main points. It is based on a Taylor expansion of the inequality $\int H_tdx\geq1$, when $t=1-s$ goes to $1$, and where $H_t$ is defined in~\eqref{eq-defHt}. Equivalently, this inequality can be written as
$$
\int t^{-n}\PAR{\inf_{h\in\R^n} \BRA{g\left(\frac{z}{t}-\frac{s}{t}h\right)+\frac{s}{t}W(h)}}^{-n}dz\geq 1. 
$$
Changing variables in the integral by letting $x=z/t$, and letting  $u=s/t$, the inequality becomes
$\int \varphi_u^{-n} \, dx \geq 1$
for any $u>0$, where  for positive $u$ $$\varphi_u(x)=\inf_{h\in\R^n} \BRA{g\big(x-u h\big)+uW(h)}.$$
 Since $\int g^{-n}dx=1$, this is
$$
\int \frac{\varphi_u^{-n}-g^{-n}}{u}dx \geq 0
$$
for any $u>0$.
The main goal is now to consider the limit as $u \to 0$, by computing the limit
\begin{equation}
\label{eq-limit}
 \lim_{u\rightarrow 0} \int \frac{\varphi_u^{-n}-g^{-n}}{u}dx. 
\end{equation}

\begin{lemma}
\label{lem-appendice}
For any $x\in\R^n$, 
$$
\lim_{u\rightarrow 0^+}\frac{\varphi_u(x)-g(x)}{u}=-W^*(\nabla g(x)).
$$
\end{lemma}
\begin{eproof}
For any $x\in\R^n$, from the definition of $\varphi_u$, we have for any $h\in\R^n$, 
$$
\frac{\varphi_u(x)-g(x)}{u}\leq \frac{g(x-uh)-g(x)}{u}+W(h)=-\nabla g(x)\cdot h+W(h)+o(u). 
$$
It follows that $\limsup_{u\rightarrow 0^+}\frac{\varphi_u(x)-g(x)}{u}\leq -\nabla g(x)\cdot h+W(h)$ for any $h$, and then by taking the infimum over $h\in\R^n$, 
$$
\limsup_{u\rightarrow 0^+}\frac{\varphi_u(x)-g(x)}{u}\leq -W^*(\nabla g(x)). 
$$
Now, one can observe that $$\varphi_u(x)=\inf_{h, \,\,uW(h)\leq g(x)} \BRA{g\big(x-u h\big)+uW(h)},$$ so that
\begin{equation}\label{eq-gphi}
\frac{g(x)-\varphi_u(x)}{u}=\!\!\sup_{h,\,\,uW(h)\leq g(x)}\!\!\BRA{\frac{g(x)-g(x-uh)}{u}-W(h)}\leq\!\!\sup_{h,\,\,uW(h)\leq g(x)}\!\!\BRA{\nabla g(x)\cdot h+|h|\varepsilon(u|h|)-W(h)},
\end{equation}
where $\varepsilon$ is an appropriate function satisfying  $\lim_{u\rightarrow0}\varepsilon(u)=0$. 

Let now $r=\sup\{u|h|,\,\,uW(h)\leq g(x)\}$. From the hypothesis (H1),  
\begin{equation}
\label{eq-ine-r}
r\leq \sup\BRA{u|h| ;\,\,\frac{u|h|^2}{2C}\leq g(x)}\leq D\sqrt{ug(x)}, 
\end{equation}
where $D$ is a constant. Generally, $D$ denotes a constant and can change from line to line. The bound~\eqref{eq-ine-r} gives 
$$
\frac{g(x)-\varphi_u(x)}{u}\leq\!\!\sup_{h,\,\,uW(h)\leq g(x)}\!\!\BRA{\nabla g(x)\cdot h+|h|\varepsilon(D\sqrt{ug(x)})-W(h)}.
$$
Let now $\eta>0$. Then there exists  $u_0>0$ such that $\forall u\in(0,u_0]$,    $\varepsilon\left(D\sqrt{ug(x)}\right)\leq \eta$, so that
$$
\frac{g(x)-\varphi_u(x)}{u}\leq\!\!\sup_{h,\,\,uW(h)\leq g(x)}\!\!\BRA{\nabla g(x)\cdot h+|h|\eta-W(h)}\leq\sup_{h\in\R^n}\BRA{\nabla g(x)\cdot h+|h|\eta-W(h)} .
$$
By (H1) the supremum is reached, say on a ball of center 0 and radius $R>0$ independent of $\eta <1$. Hence 
$$
\frac{g(x)-\varphi_u(x)}{u}\leq\sup_{h\in\R^n}\BRA{\nabla g(x)\cdot h-W(h)}+R\eta= W^*(\nabla g(x))+R\eta.
$$
The result follows by taking the superior limit and then letting $\eta$ go to 0. 
\end{eproof}

\medskip

To compute the limit~\eqref{eq-limit}, we use the dominated convergence theorem. Since everywhere  $\frac{\varphi_u^{-n}-g^{-n}}{u}$ goes to $n\frac{W^*(\nabla g)}{g^{n+1}}$ when $u\rightarrow0$, we only need to give a uniform bound (in $u$) of the quantity of $\frac{\varphi_u^{-n}-g^{-n}}{u}$. 

\medskip

For any $0<a \leq b$, the following holds $|a^{-n}-b^{-n}|\leq n|a-b|a^{-1-n}$. Since $0\leq\varphi_u(x)\leq g(x)$ by definition of $\varphi_u$,  we can apply this inequality to $a=\varphi_u(x)$ and $b=g(x)$, obtaining 
$$
\left|\frac{\varphi_u(x)^{-n}-g(x)^{-n}}{u}\right|\leq n \, \left|\frac{\varphi_u(x)-g(x)}{u}\right| \; \varphi_u(x)^{-1-n}.
$$

\medskip

{\bf Bound on $|{\varphi_u(x)-g(x)}|/{u}$:}

First, from the equality in~\eqref{eq-gphi} and a Taylor expansion, 
$$
\frac{g(x)-\varphi_u(x)}{u}\leq\!\!\sup_{h,\,\,uW(h)\leq g(x)}\!\!\BRA{\mathcal Dg(x,u|h|)|h|-W(h)},  
$$
where $\mathcal Dg(x,s)=\sup_{|x-y|\leq s}|\nabla g(y)|$. 
We assume now that $u\in]0,1]$. Then, from~\eqref{eq-ine-r},  $r\leq D\sqrt{g(x)}$. Hence, by~(H1), 
$$
\frac{g(x)-\varphi_u(x)}{u}\leq\!\!\sup_{h,\,\,uW(h)\leq g(x)}\!\!\BRA{\mathcal Dg\PAR{x,D\sqrt{g(x)}}|h|-W(h)}\leq \sup_{h \in \mathbb R^n}\BRA{\mathcal Dg\PAR{x,D\sqrt{g(x)}}|h|-\frac{|h|^2}{2C}}.  
$$
The explicit computation of the infimum gives 
$$
\frac{g(x)-\varphi_u(x)}{u}\leq D \; \mathcal Dg\PAR{x,D\sqrt{g(x)}}^2. 
$$
Then, from the hypothesis (H2), the estimation of $\mathcal Dg$ gives the bound 
$$
0\leq\frac{g(x)-\varphi_u(x)}{u}\leq D(|x|^2+1), \qquad u \in ]0, 1], x \in \mathbb R^n.
$$

{\bf Bound on $\varphi_u(x)$:}

From the hypotheses (H1) and  (H2) we have 
$$
\varphi_u(x)\geq \inf_h\BRA{\frac{1}{C}|x-uh|^2+u\frac{|h|^2}{2C}}+ \frac{1}{C}.
$$
When $u\in]0,1]$, the explicit computation of the infimum gives again
$$
\varphi_u(x)\geq D(|x|^2+1).
$$
Finally, we have obtained the upper bound 
$$
\left|\frac{\varphi_u(x)^{-n}-g(x)^{-n}}{u}\right|\leq D(|x|^2+1)^{-n}, \qquad u \in ]0, 1], x \in \mathbb R^n.
$$
The dominated convergence theorem can then be applied. The proof of Theorem~\ref{thm-super-etoile} is then complete.

 
\section{Proof of Theorem~\ref{thm-brascamp-lieb-2}}
\label{app-BL2}

We adapt the argument of~\cite{BL00}. 

We will assume throughout the proof that $V$ is $\mathcal C^3$ with bounded derivatives $\nabla^2 V$ and $\nabla^3 V$, and that there exists $\rho>0$ such that uniformly in $\R^n$, $\He{V}\geq \rho \,{\rm Id}_n$. Then the result extends to $V$ as in the Theorem by approximation. 

\smallskip

Let $f$ be a $\mathcal C^1$ compactly supported function satisfying $\int fd\mu=0$. We apply the Borell-Brascamp-Lieb inequality~\eqref{eq-borell-bl} for $t=s=1/2$,  $F=\exp(-V),$ $G=\exp(2\delta f-V)/Z_\delta$ ($\delta>0$) where $Z_\delta=\int \exp(2 \delta f)d\mu$, and finally $H=\exp(\phi_\delta-V)$ where 
\begin{equation}
\label{eq-def-phi}
\!
\phi_\delta(z)=-n\log{\inf_{h\in\R^n}\BRA{Z_\delta^{1/n}\exp\PAR{-\frac{2\delta}{n}f(z+h)+\frac{V(z+h)}{n}}+\exp\PAR{\frac{V(z-h)}{n}}}}
+n\log(2)+V(z). 
\; 
\end{equation}
Then~\eqref{eq-borell-bl} ensures that $\int  e^{\phi_\delta} d\mu\geq 1$.  The rest of the proof is devoted to a Taylor expansion of $\int  \exp(\phi_\delta)d\mu$ as $\delta $ goes to 0. 

\medskip

By convexity of $V$, for any $\delta >0$ the function in~\eqref{eq-def-phi} to be minimised is coercive, so indeed admits a (possibly non unique) minimiser, which we first estimate by giving a Taylor expansion as $\delta \to 0.$ 

For this, let $\delta >0$ be given and let $h_\delta$ be any minimiser. Then
\begin{multline}
\label{eq-sol-h}
Z_{\delta}^{1/n} \Big(-2\delta \nabla f(z+h_\delta)+\nabla V(z+h_\delta) \Big) \exp\PAR{-\frac{2\delta}{n}f(z+h_\delta)+\frac{1}{n}V(z+h_\delta)}\\=
\nabla V(z-h_\delta)\exp\PAR{\frac{1}{n}V(z-h_\delta)}.
\end{multline}

{\bf 1. In a first step} we prove that $h_\delta=O(\delta)$ uniformly in $z$: in other words, there exists a constant $C>0$ such that for any $\delta$ small enough and any $z\in\R^n$, 
$$
|h_\delta|\leq C \, \delta.
$$

For this, first, since $\int fd\mu=0$ and $f$ is compactly supported, 
\begin{equation}\label{ZZ}
Z_\delta^{1/n}=\PAR{\int e^{2 \delta f}d\mu}^{1/n}=1+\frac{2\delta^2}{n}\int f^2d\mu+o(\delta^2).
\end{equation}

Let us now assume that the support of $f$ is included in the ball $\{\vert x \vert < R\}$ with $R>0$. There are two cases, depending on whether $|z+h_\delta|\geq R$ or $|z+h_\delta| \leq R$. 

\medskip
\begin{itemize}
\item {\bf First, assume that $|z+h_\delta|\geq R$. } Then equation~\eqref{eq-sol-h} becomes 
$$
Z_{\delta}^{1/n} \nabla V(z+h_\delta) \; \exp\PAR{\frac{1}{n}V(z+h_\delta)}=
\nabla V(z-h_\delta) \; \exp\PAR{\frac{1}{n}V(z-h_\delta)}.
$$
In other words, by~\eqref{ZZ} and taking the scalar product by $h_\delta$, 
$$
\Big(\frac{2\delta^2}{n}\int f^2d\mu+o(\delta^2)\Big) \; \nabla \Phi(z+h_\delta)\cdot h_\delta=
\nabla \Phi(z-h_\delta)\cdot h_\delta -\nabla \Phi(z+h_\delta)\cdot h_\delta
$$
where $\Phi=\exp(\frac{1}{n}V)$, that is, 
$$
- \Big(\frac{2\delta^2}{n}\int f^2d\mu+o(\delta^2)\Big) \nabla V(z+h_\delta)\cdot h_\delta=
e^{-V(z + h_\delta)}  \int_{-1}^{1} h_\delta\cdot \He{\Phi}(z+th_\delta)h_\delta dt.
$$

Now $\He{V}\geq \rho  \,{\rm Id}_n$ so 
\begin{equation}\label{Hphi}
\He{\Phi} = \frac{e^V}{n} (\He{V} + \frac{1}{n} \nabla V \otimes \nabla V) \geq \frac{e^V}{n} \big(\rho \, {\rm Id}_n + \frac{1}{n} \nabla V \otimes \nabla V \big).
\end{equation} 
Hence
\begin{equation}
\label{eq-maj-h}
- \Big(2\delta^2\int f^2d\mu+o(\delta^2)\Big) \nabla V(z+h_\delta)\cdot h_\delta \geq e^{-V(z+h_\delta)}\int_{-1}^{1} \Big(\rho|h_\delta|^2+ \frac{1}{n} |\nabla V(z+th_\delta)\cdot h_\delta|^2 \Big) e^{V(z+th_\delta)} dt.
\end{equation}

In particular, $\nabla V(z+h_\delta)\cdot h_\delta \leq 0$ on the left-hand side for $\delta$ small enough, independently of $z$ since the $o(\delta^2)$ comes from $Z_{\delta}$, see~\eqref{ZZ},  and is uniform in $z$; hence for any $t\in[-1,1]$
$$
V(z+th_\delta)-V(z+h_\delta) \geq (t-1) \nabla V(z+h_\delta)\cdot h_\delta \geq 0
$$
by convexity of $V$. Moreover $\nabla V(z+th_\delta)\cdot h_\delta \leq \nabla V(z+h_\delta)\cdot h_\delta$ again by convexity, whence
$$
|\nabla V(z+th_\delta)\cdot h_\delta|\geq |\nabla V(z+h_\delta)\cdot h_\delta|.
$$

Collecting all terms,~\eqref{eq-maj-h} leads to
$$
\Big(2\delta^2 \! \int f^2d\mu+o(\delta^2)\Big) |\nabla V(z+h_\delta)\cdot h_\delta|\geq  \,2\rho |h_\delta|^2+ \frac{2}{n} |\nabla V(z+h_\delta)\cdot h_\delta|^2
\geq 4 \sqrt{\frac{\rho}{n}} |h_\delta| \,  |\nabla V(z+h_\delta)\cdot h_\delta|.
$$
for $\delta$ small enough, and where the $o(\delta^2)$ is uniform in $z$. Hence there exists a constant $A>0$~such~that 
$$
|h_\delta|\leq A \, \delta^2, 
$$
for any $\delta$ small enough and any $z$, whenever $|z+h_\delta|\geq R$. 
\medskip

\item {\bf Assume now that $|z+h_\delta|\leq R$. } 
Let us write equation~\eqref{eq-sol-h} as 
 \begin{eqnarray*}
\nabla \Phi(z+h_\delta)-\nabla \Phi(z-h_\delta)&
=&
\Big[ 1 - Z_{\delta}^{1/n} \exp\PAR{-\frac{2\delta}{n}f(z+h_\delta)} \Big] \nabla \Phi (z+h_\delta)\\
&&+\,
 2\delta Z_{\delta}^{1/n} \exp\PAR{-\frac{2\delta}{n}f(z+h_\delta)+\frac{1}{n}V(z+h_\delta)} \nabla f(z+h_\delta).
\end{eqnarray*}

Then $f$, $V$ and their gradients are continuous and then uniformly bounded on the ball $\{\vert x \vert \leq R\}$, so by~\eqref{ZZ} there exists a constant $A$ such that for all $\delta$ small enough and all $z$ with $|z+h_\delta|\leq R$
$$
 \Big|\nabla \Phi(z+h_\delta)-\nabla \Phi(z-h_\delta) \Big| \leq A \, \delta. 
$$

Hence, by the Cauchy-Schwarz inequality and the bound $\He{\Phi} \geq \frac{\rho}{n} e^V {\rm Id}_n,$ a consequence of~\eqref{Hphi},
$$
A \delta \vert h_{\delta} \vert
\geq
\Big(  \nabla \Phi(z+h_\delta)-\nabla \Phi(z-h_\delta) \Big) \cdot h_{\delta}
=
\int_{-1}^1h_{\delta} \cdot \He{\Phi}(z+th_\delta)h_\delta \, dt
\geq \frac{2 \rho}{n} e^{\min V} \vert h_{\delta} \vert^2.
$$
By uniform convexity the function $V$ is indeed bounded from below on $\mathbb R^n$, so there exists a constant $B$ such that for all $\delta$ small enough and all $z$ with $|z+h_\delta|\leq R$
$$
 |h_\delta | \leq B \, \delta.
$$
\end{itemize}

All cases being covered, our {\bf first step} is completed.

\medskip

{\bf 2. In a second step} we perform a first-order Taylor expansion of the equality~\eqref{eq-sol-h}. For $z$ fixed, it gives
\begin{equation}
\label{eq-dl1}
-\delta\nabla f(z)+\He{V}(z)h_\delta-\frac{\delta}{n}f(z)\nabla V(z)+\frac{h_\delta\cdot\nabla V(z)}{n}\nabla V(z)+o_z(\delta)=0,
\end{equation}
where $o_z(\delta)$ depends on $z$, $\delta$ and $h_\delta$. Since $ |h_\delta | \leq C \, \delta$ by the first step, uniformly in $z$, one deduces from~\eqref{eq-dl1} that 
$$
|o_z(\delta)|\leq A \; \delta^2 \; (|\He{V}(z)|+|\nabla V(z)|^2+1)
$$ 
for a constant $A$ and for any $z$.

In the sequel we let $H (z)$ denote positive polynomial functions in $V(z)$, $\nabla V(z),$ etc., independent of $\delta$ small and which can change from line to line. The latter inequality can then be written as 
\begin{equation}
\label{eq-def-oh}
|o_z(\delta)|\leq \delta^2 H(z).
\end{equation}

Let now $X=\nabla f\cdot \He{V}^{-1}\nabla V$ and $Y=\nabla V\cdot \He{V}^{-1}\nabla V.$ Taking the scalar product of~\eqref{eq-dl1} with $\He{V}^{-1}\nabla V$ one~gets 
$$
h_\delta\cdot\nabla V=\delta\frac{X+\frac{fY}{n}}{1+\frac{Y}{n}}+o_z(\delta)
$$ 
at the point $z$, where $o_z(\delta)$ satisfies~\eqref{eq-def-oh} since in particular $ \He{V}^{-1} \leq \rho^{-1}  \,{\rm Id}_n.$ Then, again by~\eqref{eq-dl1}, 
$$
h_\delta=
\delta\SBRA{\He{V}^{-1}\nabla f+\frac{\He{V}^{-1}\nabla V}{n}\frac{f-X}{1+\frac{Y}{n}}} + o_z(\delta)
$$
 where again $o_z(\delta)$ satisfies~\eqref{eq-def-oh}.

\medskip

We now compute the second-order Taylor expansion of the function $\phi_\delta$. 
First, from the expansion $e^x=1+x+x^2/2+x^3 e^{\theta x} /6$ with $\theta\in(0,1)$, we have at the point $z$, 
$$
\phi_\delta
=
-n\log (1 + \psi_{\delta})
$$
with
$$
\psi_\delta
=
-\frac{\delta}{n}f-\frac{\delta h_\delta\cdot \nabla f}{n}+\frac{h_\delta\cdot\He{V}h_\delta}{2n}+\frac{\delta^2}{n^2}f^2+\frac{(h_\delta\cdot\nabla V)^2}{2n^2}-\frac{\delta f}{n^2}h_\delta\cdot\nabla V + \frac{\delta^2}{n}\int f^2d\mu +\bar{o}_z(\delta^2).
$$
Here $\bar{o}_z(\delta^2)$ now satisfies 
\begin{equation}
\label{eq-def-ohbarre}
|\bar{o}_z(\delta^2)|\leq \delta^3 H(z)\exp\PAR{\delta^3 K_3(z)}
\end{equation}
with $|K_3(z)|\leq A(|\nabla V|+|\nabla^2 V|+|\nabla^3 V|)$ for an universal constant $A$.

\smallskip

We now observe that for small $\delta$ one has $\phi_{\delta}(z) \leq n \log 2$ for all $z$, that is, $\psi_{\delta} (z) \geq -1/2.$ Indeed, for small $\delta$ one has
$$
Z_\delta^{1/n}\exp\PAR{-\frac{2\delta}{n}f(x)} \geq \frac{1}{2}
$$
uniformly in $x \in \mathbb R^n$, by~\eqref{ZZ} and since $f$ is bounded from above. Hence for any $h \in \mathbb R^n$
\begin{multline*}
{Z_\delta^{1/n}\exp\PAR{-\frac{2\delta}{n}f(z+h)+\frac{V(z+h)}{n}}+\exp\PAR{\frac{V(z-h)}{n}}}
\\
\geq 
\frac{1}{2}\exp\PAR{\frac{V(z+h)}{n}}+\exp\PAR{\frac{V(z-h)}{n}}
\geq
\frac{1}{2} \Big[ \exp\PAR{\frac{V(z+h)}{n}}+\exp\PAR{\frac{V(z-h)}{n}} \Big]
\geq
e^{V(z)/n}
\end{multline*}
by convexity of $e^{V/n}$. The bound on $\phi_{\delta}$ follows by its definition~\eqref{eq-def-phi}.

Now from the expansion $(1+x)^{-n}=1-nx+n(n+1)x^2/2 - n(n+1)(n+2) x^3 (1+\theta x)^{-n-3}/6$ with $\theta\in (0,1)$ and~\eqref{eq-dl1}, we get  
$$
(1+\psi_{\delta})^{-n} 
= 
1  + \delta f + \frac{h_\delta\cdot\He{V}h_\delta}{2} + \delta^2 \frac{n-1}{2n}f^2+\frac{(h_\delta\cdot\nabla V)^2}{2n} - \delta^2 \! \int f^2d\mu +\bar{o}_z (\delta^2)
$$
for a $\bar{o}_z (\delta^2)$ satisfying~\eqref{eq-def-ohbarre}: here we use that $\psi_{\delta} (z) \geq -1/2$ so that $1 + \theta \psi_{\delta} \geq 1/2$ in the Taylor expansion, uniformly in $z$ and $\delta$. The above expressions of $h_\delta$ and $h_\delta\cdot\nabla V$ finally give  
$$
(1+\psi_{\delta})^{-n} 
=
1+ \delta f +\frac{\delta^2}{2}\nabla f\cdot\He{V}^{-1}\nabla f - \frac{\delta^2}{2} \frac{(f-X)^2}{n+Y} + \frac{\delta^2}{2} f^2 - \delta^2 \int f^2d\mu +\bar{o}_z (\delta^2).
$$

\medskip

In conclusion, by integration the second-order Taylor expansion of the Borell-Brascamp-Lieb inequality $\int (1+\psi_{\delta})^{-n} d\mu =  \int e^{\phi_\delta}d\mu\geq 1$ implies 
$$
\int f^2d\mu\leq \int \nabla f\cdot\He{V}^{-1}\nabla fd\mu - \int \frac{(f-X)^2}{n+Y} d\mu
$$
for all $\mathcal C^1$ compactly supported $f$ such that $\int fd\mu=0$. Here we use that $\delta^{-2} \int \bar{o}_z (\delta^2) e^{-V(z)} dz \to 0$ as $\delta \to 0$ by~\eqref{eq-def-ohbarre}, since the right-hand side in~\eqref{eq-def-ohbarre} is in $L^1(e^{-V})$ by our hypotheses on $V$. By definition of $X$ and $Y$ this concludes the argument.


\section{Link with G. Harg\'e's bound~\eqref{harge}}
\label{sec-appendix-2}

In this Appendix, we observe that G. Harg\'e's bound \eqref{harge} can be formally  recovered by linearization in the Monge-Amp\`ere equation \eqref{MA}.
 Let indeed $f$ be a smooth function such that $\displaystyle \int f \, d\mu = 0$, and $\mu_2 = (1+ \e \, f) \mu$ for $\e >0$, and expand the transport map $\nabla \varphi(x)$ sending $\mu_1 = \mu$ onto $\mu_2$ as $x + \e \nabla \theta_1 (x) + \e^2 \nabla \theta_2 (x) + o(\e^2)$. Taking logarithms in~\eqref{MA} with such $\mu_1$ and $\mu_2$ and observing that
$$
\log \det (\He{ \varphi}) = \log \det \Big(I + \e \He{\theta_1} + \e^2 \He{ \theta_2} + o(\e^2) \Big)\\
 = \e \Delta \theta_1 + \e^2 \Delta \theta_2 - \frac{\e^2}{2} \mathrm{tr} \big[( \He{ \theta_1})^2 \big] + o(\e^2),
$$
a second-order Taylor expansion ensures that $f = - L \theta_1$ in the first-order terms; moreover
$$
f^2 = - \nabla \theta_1 \cdot \He{V} \nabla \theta_1 + 2 L \theta_2 + 2 \nabla f \cdot \nabla \theta_1 - \mathrm{tr} \big[ (\He{ \theta_1})^2 \big]
$$
in the second-order terms. Assume now that $\He{V} >0$, and let $M = \He{V}^{1/2} >0$. Then
$$
- \nabla \theta_1 \cdot \He{V} \nabla \theta_1 + 2 \nabla f \cdot \nabla \theta_1 = \vert M^{-1} \nabla f \vert^2 - \vert M \nabla \theta_1 - M^{-1} \nabla f \vert^2 
$$
so that
\begin{equation}
\label{demoharge}
\int f^2 \, d\mu = \int \nabla f \cdot \He{V}^{-1} \nabla f \, d\mu - \int \Big( \vert M \nabla \theta_1 - M^{-1} \nabla f \vert^2 + \mathrm{tr} \big[ (\He{\theta_1})^2\big] \Big) \, d\mu
\end{equation}
by integration. At this point one recognizes terms in the proof of~\cite[Th.~1]{harge-jfa} : one observes that $f = - L \theta_1$ so $\nabla f = M^2 \theta_1 - X$ by differentiation, where $X \in \rr^n$ is the vector with coordinates $L(\partial_i \theta_1)$; hence
$$
\vert M \nabla \theta_1 - M^{-1} \nabla f \vert^2 = \vert M^{-1} X \vert^2 \geq \frac{1}{S} \vert X \vert^2
$$
if moreover $\He{V} \leq S$. In particular
$$
\int \vert M \nabla \theta_1 - M^{-1} \nabla f \vert^2 \, d\mu \geq \frac{1}{S} \sum_i \int \Big( L ( \partial_i \theta_1) \Big)^2 \, d\mu
 \geq \frac{R}{S} \sum_{i,j} \int \big( \partial^2_{ji} \theta_1 \big)^2 \, d\mu
$$
by \eqref{gamma2}, if $\He{V} \geq R\,{\rm Id}_n$. Hence
\begin{equation}\label{demoharge2}
\! \int \! \Big( \vert M \nabla \theta_1 - M^{-1} \nabla h \vert^2 \!+  \mathrm{tr} \big[ (\He{ \theta_1})^2\big] \! \Big)  d\mu
\geq 
\!
\left(\! 1 + \frac{R}{S} \right) \sum_{i,j} \!  \int \! \! \big( \partial^2_{ji} \theta_1 \big)^2  d\mu
\geq 
\frac{1}{n} \! \left(\! 1 + \frac{R}{S} \right) \left( \int \! \Delta \theta_1 \, d\mu \right)^2 \!
\end{equation}
since moreover by the Cauchy-Schwarz inequality
$$
\left( \int \Delta \theta_1 \, d\mu \right)^2 =  \left( \sum_i \int \partial_{ii} \theta_1 \, d\mu \right)^2 \leq n \sum_{i} \left( \int \partial_{ii} \theta_1 \, d\mu \right)^2 \leq n \sum_{i, j} \left( \int \partial_{ij} \theta_1 \, d\mu \right)^2.
$$
By \eqref{demoharge} and \eqref{demoharge2} we finally recover \eqref{harge} since by integration by parts and \eqref{gamma2}
$$
 \int \Delta \theta_1 \, d\mu = \int \nabla \theta_1 \cdot \nabla V \, e^{-V} \, dx = - \int L \theta_1 \, V \, e^{-V} \, dx = \int f \, V \, d\mu.
 $$

\end{appendix}

\bigskip

\noindent{\bf Acknowledgements.} 
{The authors are grateful to both referees for a careful reading of the manuscript and helpful comments and questions which improved the paper. In particular they raised most relevant issues. This work was partly written while the authors were visiting Institut Mittag-Leffler in Stockholm; it is a pleasure for them to thank this institution for its kind hospitality and participants for discussions on this and related works.
This research was supported  by the French ANR-12-BS01-0019 STAB project. 



\bigskip

{\footnotesize{

}


\begin{thebibliography}{10}


\bibitem{ambrosio-gigli-savare}
L.~Ambrosio, N.~Gigli and G.~Savar{\'e}.
\newblock {\em Gradient flows in metric spaces and in the space of probability
  measures}.
\newblock Lectures in Math. ETH Z\"urich. Birkh\"auser, Basel, 2008.

\bibitem{ABJ16}
M. Arnaudon, M. Bonnefont and A. Joulin.
Intertwinings and generalized Brascamp-Lieb Inequalities. 
To appear on {\em Rev. Mat. Iberoam.} 

\bibitem{BBG12}
D. Bakry, F. Bolley, and  I. Gentil.
 Dimension dependent hypercontractivity for {G}aussian kernels.
 {\em Prob. Theor. Rel. Fields}, 154(3): 845--874, 2012. 

\bibitem{bgl-book}
D.~Bakry, I.~Gentil and M.~Ledoux.
\newblock {\em Analysis and geometry of {M}arkov diffusion operators}, volume 348 of {\em Grund. Math.
  Wiss.}
\newblock Springer, Berlin, 2014.

\bibitem{bakryledoux-liyau}
D.~Bakry and M.~Ledoux.
\newblock A logarithmic {S}obolev form of the {L}i-{Y}au parabolic inequality.
\newblock {\em Rev. Mat. Iberoam.}, 22(2):683--702, 2006.

\bibitem{BK08}
F.~Barthe and A. Kolesnikov. 
Mass transport and variants of the logarithmic {S}obolev inequality. 
{\em J. Geom. Anal.}, 18(4):921--979, 2008. 


\bibitem{bgl}
S. G. Bobkov, I. Gentil and M. Ledoux.
Hypercontractivity of {H}amilton-{J}acobi equations.
{\em J. Math. Pures Appl. (9)}, 80(7):669--696, 2001.
      
\bibitem{bgrs14}      
S. Bobkov, N. Gozlan, C. Roberto and P.-M. Samson.
Bounds on the deficit in the logarithmic Sobolev inequality. 
{\em J. Funct. Anal.}, 267:4110--4138, 2014.

\bibitem{BL00}   
S. Bobkov and M. Ledoux. 
From {B}runn-{M}inkowski to {B}rascamp-{L}ieb and to logarithmic {S}obolev inequalities 
{\em Geom. Funct. Anal.}, 10:1028--1052, 2000. 

\bibitem{bobkov-ledoux-sob}   
S. Bobkov and M. Ledoux. 
From Brunn-Minkowski to sharp Sobolev inequalities.
{\em Ann. Mat. Pura Appl. (4)}, 187(3):369--384, 2008. 

\bibitem{BL09}   
S. Bobkov and M. Ledoux. 
Weighted Poincaré-type inequalities for Cauchy and other convex measures.
{\em Ann. Probab.}, 37(2):403--427, 2009. 

\bibitem{BGG11}
F.~Bolley, I.~Gentil and A.~Guillin.
\newblock Convergence to equilibrium in {W}asserstein distance for
  {F}okker-{P}lanck equations.
\newblock {\em J. Funct. Anal.}, 263(8):2430--2457, 2012.

\bibitem{BGG13}
F.~Bolley, I.~Gentil and A.~Guillin.
\newblock Dimensional contraction via Markov transportation distance. 
\newblock  {\em J. London Math. Soc.},  90(1):309--332, 2014.

\bibitem{BGGK}
F.~Bolley, I.~Gentil, A.~Guillin and K. Kuwada.
\newblock Equivalence between dimensional contractions in Wasserstein distance and curvature-dimension condition. 
To appear on {\em Annali della Scuola Norm. Sup. di Pisa.}


\bibitem{C}
E. Carlen.
Superadditivity of Fisher's information and logarithmic Sobolev inequalities.
{\em J. Funct. Anal.}, 101:194--211, 1991.


\bibitem{cordero}
D. Cordero-Erausquin.
Some applications of mass transport to {G}aussian type inequalities.
{\em Arch. Rat. Mech. Anal.}, 161:257--269, 2002. 

\bibitem{cordero15}
D. Cordero-Erausquin.
Transport inequalities for log-concave measures, quantitative forms and applications.
To appear on {\em Canadian J. Math.}

\bibitem{daneri-savare}
S. Daneri and G. Savar\'e.
Lecture notes on gradient flows and optimal transport. Optimal transportation,  
{\em London Math. Soc. Lecture Note Ser.}, 413:100--144, 2014.		

\bibitem{del-pino-dolbeault03}
M. Del Pino and J. Dolbeault. 
The optimal Euclidean $L^p$-Sobolev logarithmic inequality. 
{\em J. Funct. Anal.}, 197(1):151--161, 2003.

\bibitem{E}
R. Eldan.
A two-sided estimate for the Gaussian noise stability.
{\em Invent. Math.} 201(2):561--624, 2015.

\bibitem{EKS13}
M.~Erbar, K.~Kuwada and K.-T. Sturm.
\newblock On the equivalence of the entropic curvature-dimension condition and
  {B}ochner's inequality on metric measure spaces.
 {\em Invent. Math.} 201(3):993--1071, 2015.

\bibitem{FIL14}
M.~Fathi, E. Indrei and M. Ledoux.
Quantitative logarithmic {S}obolev inequalities and stability estimates. 
{\em Disc. Cont. Dynamical Syst.}, 36:6835--6853, 2016.
 
  \bibitem{FI13}
 A. Figalli and E. Indrei.
\newblock A sharp stability result for the relative isoperimetric inequality inside convex cones.
{\em J. Geom. Anal.}, 23:938--969, 2013.

\bibitem{FJ14}
A. Figalli and D. Jerison.
\newblock Quantitative stability for the Brunn-Minkowski inequality.
To appear on {\em Adv. Math.}

\bibitem{FMP10}
A. Figalli, F. Maggi and A. Pratelli. 
\newblock  A mass transportation approach to quantitative isoperimetric inequalities.
{\em Invent. Math.}, 182:167--211, 2010.

\bibitem{FMP13}
A. Figalli, F. Maggi and A. Pratelli. 
\newblock Sharp stability theorems for the anisotropic Sobolev and log-Sobolev inequalities on functions of bounded variation.
{\em Adv. Math.}, 242:80--101, 2013.
 
 \bibitem{gentil03}
 I. Gentil.
 The general optimal $L^p$-Euclidean logarithmic Sobolev inequality by Hamilton-Jacobi equations. 
 {\em J. Funct. Anal.} 202(2):591--599, 2003.
  
\bibitem{gentil08}
I. Gentil.
From the Pr\'ekopa-Leindler inequality to modified  logarithmic Sobolev inequality. 
{\em Ann. Fac. Sci. Toulouse,} 6, 17(2):291--308, 2008.
 
\bibitem{gentil13}
I. Gentil.
Dimensional contraction in {W}asserstein distance for diffusion semigroups  on a {R}iemannian manifold.
{\em Potential Anal.}, 42(4):861--873, 2015. 

\bibitem{gnp}
L. Goldstein, I. Nourdin and G. Peccati.
Gaussian phase transitions and conic intrinsic volumes: Steining the Steiner formula.
Preprint, 2014. 

\bibitem{harge-jfa}
G. Harg\'e.
Reinforcement of an inequality due to {B}rascamp and {L}ieb. 
{\em J. Funct. Anal.}, 254(2):267--300, 2008.

\bibitem{IM14}
E. Indrei and D. Marcon.
\newblock A quantitative log-Sobolev inequality for a two parameter family of functions.
{\em Int. Math. Res. Not.}, 20:5563--5580, 2014.

\bibitem{lisini}
S. Lisini.
Nonlinear diffusion equations with variable coefficients as gradient flows in {W}asserstein spaces.
{\em ESAIM Contr. Opt. Calc. Var.}, 15:712--740, 2009.

\bibitem{mccann97}
R. McCann.
\newblock A convexity principle for interacting gases. 
\newblock  {\em Adv. Math.}, 128(1):153--179, 1997.


\bibitem{nguyen}
V. H. Nguyen.
\newblock Dimensional variance estimates of {B}rascamp-{L}ieb type and a local approach to dimensional {P}r\'ekopa theorem.
\newblock  {\em J. Funct. Anal.}, 266(2):931--955, 2014.


\bibitem{ov00}
F.~Otto and C.~Villani.
\newblock Generalization of an inequality by {T}alagrand and links with the
  logarithmic {S}obolev inequality.
\newblock {\em J. Funct. Anal.}, 173(2):361--400, 2000.


\bibitem{talagrand96}
M.~Talagrand.
\newblock Transportation cost for {G}aussian and other product measures.
\newblock {\em Geom. Funct. Anal.}, 6(3):587--600, 1996.

\bibitem{villani-otp}
C.~Villani.
\newblock {\em Topics in Optimal transportation}, volume 58 of Grad. studies in math.
\newblock Amer. Math. Soc, Providence, 2003.

\bibitem{villani-book1}
C.~Villani.
\newblock {\em Optimal transport, Old and new}, volume 338 of {\em Grund. Math.
  Wiss.}
\newblock Springer, Berlin, 2009.


\end{thebibliography}
\end{document}